\documentclass[11pt]{article}
\usepackage[top=1.25in, bottom=1.25in, left=1.25in, right=1.25in]{geometry}
\usepackage{amssymb}
\usepackage{latexsym}
\usepackage{amsmath}
\usepackage{amsthm}
\usepackage{mathtools}
\usepackage{tikz-cd}
\usepackage{tikz}
\usepackage{stmaryrd}
\usepackage{enumitem} 
\usepackage{hyperref}
\usepackage{graphicx}
\usepackage{caption}
\usepackage{multicol}
\usepackage{subcaption}
\usepackage[title]{appendix}
\usepackage{mathrsfs}
\usepackage{bbm}
\usepackage{bm}
\usepackage[T1]{fontenc}
\usepackage{color,soul}
\usepackage{stackengine, scalerel}
\usepackage[font=footnotesize]{caption}

\DeclareMathOperator{\im}{im}

\graphicspath{ {./images/} }
\hypersetup{
    colorlinks,
    citecolor=blue,
    filecolor=red,
    linkcolor=blue,
    urlcolor=blue
}

\newtheorem{proposition}{Proposition}[section]

\theoremstyle{definition}
\newtheorem{definition}{Definition}[section]

\theoremstyle{definition}
\newtheorem{example}{Example}[section]

\newcommand{\real}{\mathbb R}

\newcommand{\intg}{\mathbb Z}

\newcommand{\lgl}{\langle}
\newcommand{\rgl}{\rangle}

\title{Persistent sheaf Laplacians}
\author{Xiaoqi Wei$^{1}\footnote{
 E-mail: weixiaoq@msu.edu} $, 
 and Guo-Wei Wei$^{1,2,3}$\footnote{
 E-mail: weig@msu.edu} \\% Author name
$^1$ Department of Mathematics, \\
Michigan State University, MI 48824, USA.\\
$^2$ Department of Electrical and Computer Engineering,\\
Michigan State University, MI 48824, USA. \\
$^3$ Department of Biochemistry and Molecular Biology,\\
Michigan State University, MI 48824, USA. \\
}
\date{}

\begin{document}
\maketitle
\begin{abstract}
   Recently various types of topological Laplacians have been studied from the perspective of data analysis. The spectral theory of these Laplacians has significantly extended the scope of algebraic topology and data analysis. Inspired by the theory of persistent Laplacians and cellular sheaves, this work develops the theory of persistent sheaf Laplacians for cellular sheaves, and describes how to construct sheaves for a point cloud where each point is associated with a quantity that can be devised to embed physical properties. As a result, the spectra of persistent sheaf Laplacians encode both geometrical and non-geometrical information of the given point cloud. The theory of persistent sheaf Laplacians is an elegant method for fusing different types of data and has huge potential for future development.    
\end{abstract}


\newpage
 
\tableofcontents

\newpage

\section{Introduction}
Recent years have witnessed a dramatic growth of research interest in topological data analysis (TDA) \cite{kaczynski2004computational,wasserman2018topological}, 
driven by its success in science and technology, particularly in computational biology and computer-aided drug design, one of the most challenging scientific fronts of the 21$^{st}$ century \cite{nguyen2019mathematical,nguyen2020mathdl}. 
The main workhouse of TDA is persistent homology \cite{cerri2013betti, edelsbrunner2008persistent, zomorodian2005computing}, a new branch of algebraic topology 
that is able to capture multiscale topological information of data. 
Topological deep learning was proposed  in 2017 to integrate  persistent homology and artificial intelligence (AI)   \cite{cang2017topologynet}. It has had  enormous success in data science  \cite{ meng2020weighted, nguyen2020review,townsend2020representation}.  

Persistent homology has many limitations. 
%For example, it cannot capture the shape evolution of data that does not change the topological invariant.  
For example, it cannot capture the shape evolution of data that does not change topologically.
Evolutionary de Rham-Hodge theory, or persistent Hodge Laplacian, was proposed in 2019 to overcome this limitation \cite{chen2019evolutionary}.   
The de Rham-Hodge theory connects algebraic topology with differentiable geometry, partial differential equations, harmonic analysis, and algebraic geometry \cite{dodziuk1977rham,zhao2020rham}. 
Persistent Hodge Laplacian defined on a family of differentiable manifolds offers a powerful multiscale spectral analysis of volumetric data \cite{chen2019evolutionary}.  
Additionally, persistent spectral graph theory, also called persistent combinatorial Laplacians (or simply persistent Laplacians), 
was also introduced in 2019 to overcome the drawback of  persistent homology \cite{wang2020persistent}. 
{ Combinatorial Laplacians enable us to perform spectral analysis for simplicial complexes constructed from a point cloud \cite{chung1997spectral, grone1990laplacian, lim2020hodge}.}  
Persistent combinatorial Laplacians have stimulated much interest in the past few years \cite{memoli2022persistent, wang2021hermes, liu2023algebraic} and has had success in various applications \cite{meng2021persistent,chen2022persistent,qiu2023persistent}. Both persistent Hodge  Laplacians and persistent combinatorial Laplacians are  persistent topological Laplacians, an emerging research topic in TDA.  
However, these PTLs cannot deal with labeled data, or the local non-geometric properties at individual data points.  

% aim and motivation
% I don't see very often 'mathematically' in the beginning of a sentence
 {A natural extension of the aforementioned algebraic topology, differential geometry, and graph theory approaches for data science is the sheaf theory. }
A sheaf is a mathematical tool for systematically tracking data (or states), or extracting spatial-temporal patterns, or organizing data (states), or fusing data (states) according to certain physical/mathematical rules or structures. 
It defines a relationship for adjacent data points lying in a topological space, which may be abstracted from point clouds, manifolds, or spatial-temporal data. 
Such relationships become particularly valuable when they are devised to embed certain physical laws governing the underlying data. 
A particular class of sheaves, namely cellular sheaves, has attracted much attention in the past decade for its application potentials. 
The notion of a \emph{cellular sheaf} was first introduced by Shepard in his Ph.D. thesis in 1985 \cite{shepard1985cellular}. 
In the past decade, many researchers have revived and expanded the cellular sheaf theory with applications in science and engineering \cite{curry2014sheaves, hansen2019toward, topologicalsignalprocessing}. 
Tools from sheaf theory have also been developed to study persistence modules \cite{berkouk2021derived, Bubenik_2021, curry2014sheaves, hang2021correspondence}. 
Roughly speaking, a cellular sheaf consists of a simplicial complex, an assignment of vector spaces for each simplex, and a definition of linear maps for each face relation, satisfying certain rules so that it gives rise to a \emph{sheaf cochain complex}. 
Therefore, one can define sheaf cohomology, which reflects the property of the sheaf. 
Spectral sheaf theory extends spectral graph theory to cellular sheaves, leading to sheaf Laplacians \cite{hansen2019toward}. 

The aim of this paper is to introduce { persistent sheaf Laplacians} (PSL) as an extension of persistent Laplacians to the setting of cellular sheaves.   
This work is also motivated by the need to elegantly fuse geometric and non-geometric information as in the persistent cohomology  \cite{cang2020persistent}. 
In this approach, the geometrical shape can be successfully captured by persistent homology and particularly, by persistent Laplacians, while local non-geometric properties, such as the atomic biophysical properties of a molecule and the degree of a vertex in a network,  requires additional treatment. 
It is beneficial to have mathematical tools that can distinguish between data that have very similar geometrical shapes but carrying different non-geometrical information.
From a mathematical perspective, the extension of persistent Laplacian is natural, 
since the key theorem in the theory of persistent Laplacians is indeed true for more general (co)chain complexes. 
The rest of this paper is organized as follows. In Section \ref{preliminaries}, we introduce the basics of cellular sheaf theory and discuss how to define a sheaf on a labeled simplicial complex (i.e., a simplicial complex with a quantity associated to each vertex). 
In Section \ref{persistent sheaf laplacian}, we define the $q$-th persistent sheaf Laplacian for a pair of complexes $X,Y$. 
In Section \ref{experiments}  we demonstrate the spectra of persistent sheaf Laplacians for several examples of applications.

\section{Preliminaries}\label{preliminaries}

In this section, we give a very concise introduction to cellular sheaves and persistent Laplacians.
We assume that the reader is familiar with simplicial homology.
{ We follow the standard notational convention to distinguish between chain complexes and cochain complexes.}
Most researchers define cellular sheaves on regular cell complexes \cite{curry2014sheaves,  elementaryappliedtopology, hansen2019toward, topologicalsignalprocessing}. 
For the sake of simplicity, we only discuss celluar sheaves on simplical complexes, as they are indeed regular cell complexes.

\subsection{Cellular sheaves}

%We assume the reader is familiar with simplicial homology. 
%A good reference of simplicial homology is 
%Main references of this section are \cite{JustinCurry,  elementaryappliedtopology, topologicalsignalprocessing}.
\begin{definition}
    %\cite[Chapter 9]{elementaryappliedtopology}
    A \emph{cellular sheaf} $\mathscr{S}$ on a simplicial complex $X$ consists of the following data:
 
    (1) a simplicial complex $X$, where the face relation that
    $\sigma$ is a face of $\tau$ is denoted by $\sigma \leqslant \tau$, and 

    (2) an assignment to each simplex $\sigma$ of $X$ a (finite dimensional) vector space $\mathscr{S}(\sigma)$
    and to each face relation $\sigma \leqslant \tau$ a linear morphism of vector spaces denoted by
    $\mathscr{S}_{\sigma \leqslant \tau}$
    or $\mathscr{S}(\sigma \leqslant \tau): \mathscr{S}(\sigma) \to \mathscr{S}(\tau)$, satisfying the rule 
    \begin{align*}
        \rho \leqslant \sigma \leqslant \tau \Rightarrow \mathscr{S}_{\rho \leqslant \tau} =   \mathscr{S}_{\sigma \leqslant \tau} \circ \mathscr{S}_{\rho \leqslant \sigma}
    \end{align*}
    and $\mathscr{S}_{\sigma \leqslant \sigma}={\rm id}$ is the identity map. 

    The vector space $\mathscr{S}(\sigma)$ is refereed to as the \emph{stalk} of $\mathscr{S}$ over $\sigma$, and the linear morphism 
    $\mathscr{S}_{\sigma \leqslant \tau}$ is referred to as the \emph{restriction map} of the face relation $\sigma \leqslant \tau$. 
\end{definition}

\begin{example}
    Let $X$ be a finite simplicial complex. 
    We attach to every simplex of $X$ a fixed vector space $V$ and let every restriction map be the identity map. 
    This sheaf is referred to as the \emph{constant sheaf} $\underline{V}$ on $X$.
\end{example}

\begin{definition}
    Suppose that $f: X \to Y$ is a simplicial map \cite{Munkres1984} and that $\mathscr{S}$ is a cellular sheaf on Y. The \emph{pullback sheaf} $f^{\ast} \mathscr{S}$ on X is
    given by 
    \begin{align*}
        (f^{\ast} \mathscr{S})(\sigma) = \mathscr{S}(f(\sigma)),
    \end{align*}
    and for the face relation $\sigma \leqslant \tau$ of $X$,
    \begin{align*}
        (f^{\ast}\mathscr{S})_{\sigma \leqslant \tau} = \mathscr{S}_{f(\sigma) \leqslant f(\tau)}. 
    \end{align*}
\end{definition}

\begin{example}
    Suppose that $X$ is a subcomplex of $Y$, and $\mathscr{S}$ is a sheaf on $Y$. We can define a sheaf $\mathscr{T}$ on $X$ using the data of $Y$.
    For $\sigma \in X$, let $\mathscr{T}(\sigma) = \mathscr{S}(\sigma)$. For the face relation $\sigma \leqslant \tau$ in $X$, 
    let $\mathscr{T}(\sigma \leqslant \tau) = \mathscr{S}(\sigma \leqslant \tau)$. The sheaf $\mathscr{T}$ is a pullback of $\mathscr{S}$.
\end{example}

\begin{definition}
    Suppose that $\mathscr{S}$ is a sheaf.
    A \emph{global section} $s$ of $\mathscr{S}$ is an assignment to each simplex $\sigma$ an element $s_{\sigma} \in \mathscr{S}(\sigma)$ such that 
    $\mathscr{S}_{\sigma \leqslant \tau} (s_{\sigma}) = s_{\tau}$ for any face relation $\sigma \leqslant \tau$. 
    The set of global sections is denoted by $\Gamma(X; \mathscr{S})$.
\end{definition}

\begin{example}
    This example is related to quantum physics.
    The general form of the Schr\"odinger equation for an isolated quantum system is
    \begin{align*}
        i\hbar \frac{\partial}{\partial t} | \psi(t) \rangle = H | \psi(t) \rangle
    \end{align*}
    where $t$ is time, $| \psi(t) \rangle=| \psi({\bf r}, t) \rangle$ is the state vector of the quantum system that belongs to a Hilbert space $\mathcal{H}$, ${\bf r}$ is the position vector,     and $H$ is the Hamiltonian consisting of kinetic energy and potential energy operators. In the position space representation, the kinetic energy operator is given by the Laplacian. 
    In the Schr\"odinger representation, $H$ is independent of time ($\frac{\partial H}{\partial t}=0$), and the time evolution of $| \psi(t) \rangle$ is given by 
    \begin{align*}
        | \psi(t) \rangle = e^{-iH(t-t_0)/\hbar} | \psi(t_0) \rangle,
    \end{align*}
	where $e^{-iH(t-t_0)/\hbar} $ is known as the time-evolution operator. 
    %When $\mathcal H$ is finite dimensional, 
	We can define a cellular sheaf $\mathscr{S}$   as follows.
    First $\mathbb R$ is seen as a simplicial complex such that its vertices are integers and its edges are intervals $[n,n+1]$ where $n\in \mathbb Z$.
    We attach $\mathcal H$ to each simplex. Restriction maps $\mathscr{S}_{n \leqslant [n, n+1]}$ and $\mathscr{S}_{n+1 \leqslant [n,n+1]}$
    are defined by linear maps $e^{\frac{-iH }{2 \hbar}}$ and $e^{\frac{iH }{2 \hbar}}$.
    The assignment $n \to |\psi(n)\rangle$ for all $n$ is a global section of sheaf $\mathscr{S}$.   
	%The time-reversed propagation of the quantum system can be also defined similarly. 
	Similarly, cellular sheaves can be defined for many other linear (partial) differential equations. 
    Note that in general there is no need to work with $\intg$. 		
\end{example}

\subsection{Sheaf cohomology}

For a sheaf $\mathscr{S}$ on a finite simplicial complex $X$, 
we can construct the \emph{sheaf cochain complex} of $\mathscr{S}$ as follows.
Let the $q$-th cochain group $C^q(X; \mathscr{S})$ be the direct sum of
$\mathscr{S}(\sigma)$ over all $q$-simplices $\sigma$. 
%(sometimes we omit the $\mathscr{S}$ in $C^q(X; \mathscr{S})$ and only write $C^q(X)$). 
To define the coboundary map $d$, we need a \emph{signed incidence relation} \cite{curry2014sheaves}.

\begin{definition}
    A \emph{signed incidence relation} is an assignment to every face relation $\sigma \leqslant \tau$ an integer 
    $[\sigma: \tau]$ satisfying the following conditions:

    (1) if $\dim \tau - \dim \sigma > 1$, then $[\sigma:\tau] = 0$; and

    (2) if $\gamma \leqslant \tau$ and $\dim \tau - \dim \gamma=2$, the sum $\sum_{\sigma} [\gamma: \sigma][\sigma: \tau] = 0$.
\end{definition}

If a signed incidence relation is given,  
we then define the coboundary map $d_{q}: C^q(X; \mathscr{S}) \to C^{q+1}(X; \mathscr{S})$ by 
\begin{align*}
    d^q \vert_{\mathscr{S}(\sigma)} = \sum_{\sigma \leqslant \tau} [\sigma: \tau] \mathscr{S}_{\sigma \leqslant \tau}.
\end{align*} 
Since $d^q$ is a linear morphism, its action on each stalk $\mathscr{S}({\sigma})$ determines itself.
We can verify that $d^qd^{q-1}=0$ \cite[Lemma 6.2.2]{curry2014sheaves}, so there is the sheaf cochain complex
\[
    \begin{tikzcd}
        \centering
    0 \arrow[r] & 
    C^0(X; \mathscr{S}) \arrow[r, "d"] & 
    C^1(X; \mathscr{S}) \arrow[r, "d"] & 
    C^2(X; \mathscr{S}) \arrow[r, "d"] & 
    \cdots.
    \end{tikzcd} 
\]
The $q$-th sheaf cohomology group $H^q(X; \mathscr{S})$ is defined by $\ker d^q/\im d^{q-1}$. 

A natural signed incidence relation exists for every oriented simplicial complex. 
Recall that the orientation of a simplex is determined by the ordering of its vertices. 
For an oriented simplex $\tau=[v_0, v_1, \dots, v_n]$ and its oriented face 
$\sigma =[v_0, \dots, \hat{v}_{i}, \dots, v_n]$, we let $[\sigma:\tau] = (-1)^i$.
If $\sigma$ or $\tau$ is oriented alternatively,  we let $[\sigma:\tau] = (-1)^{i+1}$.
This signed incidence relation is used throughout this paper. 
In practice we can orient a simplicial complex by a global ordering of vertices.
We remind the reader that we do not need orientation information to define a sheaf. 

\begin{figure}[htbp] 
    \centering
    \begin{subfigure}{0.5\textwidth}
        \centering
        \includegraphics[width=0.6\textwidth]{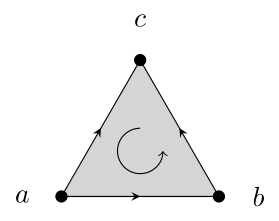}
        %\caption{Pentagon}
    \end{subfigure}\hfill
    \begin{subfigure}{0.5\textwidth}
        \centering
        \includegraphics[width=0.6\textwidth]{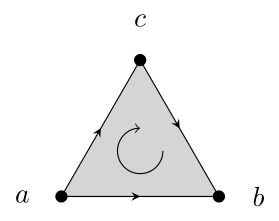} 
        %\caption{Heptagon}
    \end{subfigure}\hfill
    \caption{Different orderings determine different orientations.
    On the left, the ordering is $a,b,c$, and the oriented simplices are $[a, b], [b,c],[a,c],[a,b,c]$.
    On the right, the ordering is $a,c,b$, and the oriented simplices are $[a, b], [c, b],[a,c], [a,c,b]$.
    In practice, one would label vertices by natural numbers, and orient simplices such that labels are increasing.
    For example, if we label $a,b,c$ by $0,1,2$, then we will get the left orientation. For the sake of simplicity, we often denote a vertex by its index, so $v_i$ is written as $i$;
    for edges, we just write $ij$ instead of $[v_i, v_j]$, and adopt similar notations for higher dimensional simplices.}
    \label{orientation}
\end{figure}

\begin{example}
    We examine the sheaf cochain complex of the constant sheaf $\underline{\mathbb{R}}$ on a finite simplicial complex $X$.
    The dimension of the $q$-th cochain group $C^q(X; \mathscr{S})$ is equal to the number of $q$-simplices. 
    We use simplices to distinguish different stalks $\real$ over different simplices. 
    The $q$-th sheaf cochain group has the canonical basis $\{ \sigma \vert \dim \sigma = q \}$,
    and
    \begin{align*}
        d^q \sigma = \sum_{\sigma \leqslant \tau} [\sigma: \tau]\tau.
    \end{align*}
    %We see that in this case the sheaf cohomology coincides with the cellular cohomology of $X$.
    Suppose $X$ is $\{0,1,2,01,02,12,012\}$,
    \begin{align*}
        \begin{tikzpicture}
            \centering
            \draw[thick] (0,0) -- (2,0) -- (1, 1.732) -- (0,0);
            \draw[fill=black!20!white] (0,0) -- (2,0) -- (1, 1.732) -- (0,0);
            \node at (-0.5,0) {$0$};
            \node at (2.5,0) {$1$};
            \node at (1,2.232) {$2$};
            \draw[fill=black] (0,0) circle (2pt);
            \draw[fill=black] (1,1.732) circle (2pt);
            \draw[fill=black] (2,0) circle (2pt);
            %\caption{}
        \end{tikzpicture}
    \end{align*} 
    then the matrix representation of $d^0$ is 
    \begin{align*}
        \bordermatrix{
            ~  &   0  & 1   & 2  \cr
            01 & -1 & 1 & 0 \cr
            12 & 0  & -1 & 1 \cr
            02 & -1 & 0 & 1 \cr
        }.
    \end{align*}
    and the matrix representation of $d^1$ is 
    \begin{align*}
        \bordermatrix{
            ~  &   01  & 02   & 12  \cr
            012 & 1 & -1 & 1 \cr
        }.
    \end{align*}
    In fact for any finite simplicial complex $X$, 
    the sheaf cochain complex of $\underline{\mathbb{R}}$ coincides with the dual of the simplicial chain complex of $X$ with coefficient ring $\mathbb R$.
\end{example}

The following fact is well-known and we omit the proof.
\begin{proposition}\cite[Lemma 9.5]{elementaryappliedtopology}
    $H^{0}(X; \mathscr{S}) = \Gamma(X; \mathscr{S})$.
\end{proposition}

\subsection{Cellular sheaves on a labeled simplicial complex} \label{cellular sheaves on a labeled simplicial complex}

Suppose that there is a $1$-dimensional simplicial complex (i.e., graph) $X$ where
each vertex $v_i$ is associated with a quantity $q_i\in \real$.
Denote the edge connecting $v_i$ and $v_j$ by $e_{ij}$.
We can define a sheaf $\mathscr{S}$ on $X$ such that each stalk is $\mathbb R$, and
for the face relation $v_i \leqslant e_{ij}$, 
the morphism $\mathscr{S}_{v_i \leqslant e_{ij}}$ is the multiplication by $q_j/r_{ij}$ 
where $r_{ij}$ is the length of $e_{ij}$.
This sheaf is inspired by molecular biology and chemistry. Given a molecule, each atom and its atomic partial charge can be seen as a vertex $v_i$ and 
the associated quantity $q_i$, respectively. Conceptually, one may build a $1$-dimensional Rips or alpha complex from the point cloud of vertices or just take chemical bonds as edges.  
The assignment $q_i \to v_i$ and $q_iq_j/r_{ij} \to e_{ij}$ is a global section, 
since $\mathscr{S}_{v_i \leqslant e_{ij}}(q_i) = \mathscr{S}_{v_j \leqslant e_{ij}}(q_j) = q_iq_j/r_{ij}$. 
The quantity $q_iq_j/r_{ij}$ is the potential energy. 

The above sheaf can be generalized to high dimensional simplicial complexes (cf. \cite{wu2018weighted}).
Suppose we have the following data: a simplicial complex $X$, a set of nonzero $q_i \in \mathbb{R}$ 
associated to each vertex $v_i$, and a nowhere zero function $F: X \to \mathbb{R}$. 
We can define a sheaf where each stalk is $\mathbb{R}$, and
for the face relation $[v_0, \dots, v_n] \leqslant [v_0, \dots, v_n, v_{n+1} \dots, v_m]$,
the linear morphism $\mathscr{S}([v_0, \dots, v_n] \leqslant [v_0, \dots, v_n, v_{n+1} \dots, v_m])$ is 
the scalar multiplication by
\begin{align*}
     \frac{F([v_0, \dots, v_n])q_{n+1}\cdots q_{m}}{F([v_0, \dots, v_n, v_{n+1}, \dots, v_m])}.
\end{align*}
This is indeed a sheaf since if we have $[v_0, \dots, v_n] \leqslant [v_0, \dots, v_m] \leqslant [v_0, \dots, v_l]$, then 
\begin{align*}
    \frac{F([v_0, \dots, v_m])q_{m+1}\cdots q_l}{F([v_0, \dots, v_l])} \frac{F([v_0, \dots, v_n])q_{n+1}\cdots q_{m}}{F([v_0, \dots, v_m])} = \frac{F([v_0, \dots, v_n])q_{n+1}\dots q_l}{F([v_0, \dots, v_l])}.
\end{align*}
The assignment $q_{i_0}\cdots q_{i_n}/F([v_{i_0}, \dots, v_{i_n}]) \to [v_{i_0}, \dots, v_{i_n}]$ is a nontrivial global section.

\begin{example} \label{example2.5}
    Consider the oriented simplicial complex $\{v_0,v_1,v_2,v_0v_1,v_1v_2,v_0v_2,v_0v_1v_2\}$
    \[
        \begin{tikzpicture}
            \centering
            \draw[thick] (0,0) -- (2,0) -- (1, 1.732) -- (0,0);
            \draw[fill=black!20!white] (0,0) -- (2,0) -- (1, 1.732) -- (0,0);
            \node at (-1,0) {$(v_0, q_0)$};
            \node at (3,0) {$(v_1, q_1)$};
            \node at (1,2.232) {$(v_2, q_2)$};
            \draw[fill=black] (0,0) circle (2pt);
            \draw[fill=black] (1,1.732) circle (2pt);
            \draw[fill=black] (2,0) circle (2pt);
            %\caption{}
        \end{tikzpicture}
    \] 
    where $q_i\in \real$ is associated to $v_i$. Let $r_{01}, r_{12}, r_{02}$ be the lengths of $e_{01}, e_{12}, e_{02}$. 
    We can define the above sheaf on this complex where $F$ maps every vertex to 1, every edge $e_{ij}$ to its length $r_{ij}$, 
    and the 2-simplex $[v_0,v_1,v_2]$ to $r_{01}r_{12}r_{02}$. 
    The matrix representation of $d^0$ is
    \begin{align*}
        \bordermatrix{
            ~  &   v_0  & v_1   & v_2  \cr
            v_0v_1 & -q_1/r_{01} & q_0/r_{01} & 0 \cr
            v_1v_2 & 0  & -q_2/r_{12} & q_1/r_{12} \cr
            v_0v_2 & -q_2/r_{02} & 0 & q_0/r_{02} \cr
        },
    \end{align*} 
    and the matrix representation of $d^1$ is 
    \begin{align*}
        \bordermatrix{
            ~               &   v_0v_1  & v_0v_2 & v_1v_2  \cr
            v_0v_1v_2 & \frac{q_2}{r_{02}r_{12}} & \frac{-q_1}{r_{01}r_{12}} & \frac{q_0}{r_{01}r_{02}} \cr
        }.
    \end{align*}    
\end{example}
Note that many alternative sheaf constructions are available by appropriate choices of  $F$. For example, $F$ may map a 2-cell to the sum of its edge lengths. In practical applications, a good choice of $F$ can embed physical information into spectral representation.

\subsection{Sheaf Laplacian}

If cochain groups of a cochain complex
\[
    \begin{tikzcd}[column sep = large]
        \centering
    \cdots \arrow[r, "d^{q-2}"] & 
    A^{q-1} \arrow[r, "d^{q-1}"] & 
    A^q \arrow[r, "d^q"] & 
    A^{q+1} \arrow[r, "d^{q+1}"] & 
    \cdots
    \end{tikzcd} 
\]
are all finite dimensional inner product spaces,
the $q$-th combinatorial Laplacian $\Delta_q: A^{q} \to A^{q}$ is defined by
\begin{align*}
    \Delta_q = (d^{q})^{\ast}d^q + d^{q-1}(d^{q-1})^{\ast},
\end{align*}
where $(d^q)^{\ast}$ is the adjoint of $d^q$, and it is well-known that the kernel of $\Delta_q$ is isomorphic to the $q$-th cohomology group $H^q$.
Hansen and Ghrist \cite{hansen2019toward} applied this construction to sheaf cochain complexes and the resulting new combinatorial Laplacian is referred to as the sheaf Laplacian. 
If every stalk of a sheaf $\mathscr{S}$ is a finite dimensional inner product space, 
we can equip an inner product structure on every $C^q(X; \mathscr{S})$
such that $\mathscr{S}(\sigma)$ and $\mathscr{S}(\sigma')$ are orthogonal if $\sigma \neq \sigma'$.

\begin{example}
    We consider the $F$ defined in the same way as in Example \ref{example2.5} to evaluate the spectra of sheaf Laplacians.
    Consider the 2-simplex
    \[
        \begin{tikzpicture}
            \centering
            \draw[thick] (0,0) -- (2,0) -- (1, 1.732) -- (0,0);
            \draw[fill=black!20!white] (0,0) -- (2,0) -- (1, 1.732) -- (0,0);
            \node at (-1,0) {$(v_0, q_0)$};
            \node at (3,0) {$(v_1, q_1)$};
            \node at (1,2.232) {$(v_2, q_2)$};
            \draw[fill=black] (0,0) circle (2pt);
            \draw[fill=black] (1,1.732) circle (2pt);
            \draw[fill=black] (2,0) circle (2pt);
        \end{tikzpicture}
    \] 
    whose edges are all of length 1. The $\Delta_0$ is 
    \begin{align*}
        \begin{pmatrix}
            q_1^2+q_2^2 & -q_0q_1 & -q_0q_2 \\
            -q_0q_1  & q_0^2+q_2^2  & -q_1q_2 \\
            -q_0q_2 & -q_1q_2 & q_0^2+q_1^2
        \end{pmatrix}        
    \end{align*}
    and its eigenvalues are $\{ q_0^2+q_1^2+q_2^2, q_0^2+q_1^2+q_2^2, 0 \}$  and  the corresponding eigenvectors are
    $(-q_1/q_0, 1, 0)^T, (-q_2/q_0, 0, 1)^T$, and $(q_0/q_2, q_1/q_2, 1)^T$.   
Moreover,     $\Delta_1$ is 
    \begin{align*}
        \begin{pmatrix}
            q_0^2+q_1^2+q_2^2 & 0 & 0 \\
            0 & q_0^2+q_1^2+q_2^2 & 0 \\
            0 & 0 & q_0^2+q_1^2+q_2^2
        \end{pmatrix}        
    \end{align*}
    and its only eigenvalue is $q_0^2+q_1^2+q_2^2$. 

This example shows that the eigenvalues of 	$\Delta_0$ and 	$\Delta_1$ are  dependent on the amplitude of $q_i$, which allows the embedding of non-geometric information  in practical applications. 
However, they are not sensitive to the sign of $q_i$. 
Therefore, a (persistent) sheaf Dirac theory as an extension of recent Dirac formulation or quantum persistent homology \cite{ameneyro2022quantum} may enable us to further eliminate the sign degeneracy.      
		
Furthermore, the eigenvectors of $\Delta_0$ depend on the signs of  of $q_i$. 
As such, { it is advantageous to use sheaf eigenvectors to embed physical properties of a dataset and utilize them in machine learning. }

\end{example}

\section{Persistent Sheaf Laplacians} \label{persistent sheaf laplacian}

We review persistent combinatorial Laplacians before introducing persistent sheaf Laplacians. 

\subsection{Persistent Laplacians}

Given two chain complexes $(V,d)$ and $(W,d)$, if $V_q$, as a inner product space, is a subspace of the finite dimensional inner product space $W_q$ for all $q$, 
and the boundary operator of $(V,d)$ inherits from $(W,d)$,
then we have the following commutative diagram
\[
    \begin{tikzcd}[column sep = large]
        \cdots \arrow[r, "d_{q+2}"]
        & V_{q+1} \arrow[r, "d_{q+1}"] \arrow[d, hook, dashed] 
          & V_{q} \arrow[r, "d_q"] \arrow[d, hook, dashed] 
            & V_{q-1} \arrow[r, "d_{q-1}"] \arrow[d, dashed, hook]
              & \cdots
            \\
        \cdots \arrow[r, "d_{q+2}"]
        & W_{q+1} \arrow[r, "d_{q+1}"] 
          & W_q \arrow[r, "d_q"] 
            & W_{q-1} \arrow[r, "d_{q-1}"]
              & \cdots 
    \end{tikzcd}
\]
where dashed hooked arrows represent the inclusion map $\iota: V \hookrightarrow W$.
If we use superscripts $V,W$ to distinguish among maps and subspaces of $(V,d)$ and $(W,d)$, 
then the $q$-th persistent homology group $H_q(V,W)$ is defined by 
\begin{align*}
    \frac{Z^V_q}{B_q^W \cap Z_q^V} \cong \iota^{\bullet} (H_q(V)),
\end{align*}
whose dimension is called the $q$-th persistent Betti numbers of the pair $V,W$.
Let $\Theta_{q+1}^{V,W} = \{x \in W_{q+1} \mid d_{q+1}^Wx \in V_q \}$. 
In other words, $\Theta_{q+1}^{V,W}$ consists of chains of $W_{q+1}$ whose boundaries are in $V_q$ (actually in $Z^V_q$, since $d^2=0$),
and we have $\im \Theta_{q+1}^{V,W}=B_q^W \cap Z_q^V$. 
%If $W_q$ is an inner product space for all $q$, and $V_q$ inherits its inner product from $W_q$,
If we denote the adjoint of the inclusion $V\hookrightarrow W$ by $\pi$ (which is a projection), 
the $q$-th persistent Laplacian $\Delta_q^{V,W}$ \cite{liu2023algebraic, memoli2022persistent, wang2020persistent} is defined by 
\begin{align*}
    \Delta_q^{V,W} = (d_{q}^V)^{\ast}d_q^V + \pi d_{q+1}^W\vert_{\Theta_{q+1}^{V,W}} (\pi d_{q+1}^W\vert_{\Theta_{q+1}^{V,W}})^{\ast}.
\end{align*}
Since we know that $\ker \Delta_q^{V,W} \cong H_q(V,W)$ \cite{liu2023algebraic, memoli2022persistent, wang2020persistent},
the theory of persistent Laplacian offers an alternative method of computing persistent Betti numbers. 
Besides that, non-zero eigenvalues and eigenvectors of a persistent Laplacian contain extra information that can to be utilized.
%and encode extra information of in non-zero eigenvalues and eigenvectors.

\subsection{Persistent sheaf cohomology and persistent sheaf Laplacians}

Persistent sheaf cohomology is known to many researchers \cite{michaelrobinsontutorial, russold2022persistent, yegnesh2016persistence}.
Given two oriented simplicial complexes $X,Y$, if $X \subset Y$ and the orientation of $X$ is identical to $Y$, 
let sheaf $\mathscr{F}$ on $X$ be the pullback of the sheaf $\mathscr{G}$ on $Y$,
then we have the following commutative diagram 
\[
    \begin{tikzcd}
        \cdots \arrow[r, "d"]
         & C^{q-1}(X; \mathscr{F}) \arrow[r, "d"] 
          & C^q(X; \mathscr{F}) \arrow[r, "d"] 
           & C^{q+1}(X; \mathscr{F}) \arrow[r, "d"]
            & \cdots\\
        \cdots \arrow[r, "d"]
         & C^{q-1}(Y; \mathscr{G}) \arrow[r, "d"] \arrow[u, "\pi"]  
          & C^q(Y; \mathscr{G})    \arrow[r, "d"] \arrow[u, "\pi"]  
           & C^{q+1}(Y; \mathscr{G}) \arrow[r, "d"] \arrow[u, "\pi"]
            & \cdots 
    \end{tikzcd}
\]
where $\pi: C^{q}(Y; \mathscr{G}) \to C^q(X; \mathscr{F})$ is a projection map such that
$\pi\mid_{\mathscr{G}(\sigma)}$ is the identity map if $\sigma \in X$, and $\pi\mid_{\mathscr{G}(\sigma)} = 0$ otherwise.
Since $\pi$ is a cochain map, it induces a map $\pi^{\bullet}$ between sheaf cohomology groups of $\mathscr{F}$ and $\mathscr{G}$,
and the $q$-th persistent sheaf homology group is defined by
\begin{align*}
    \pi^{\bullet} (H^q(Y; \mathscr{G}))
\end{align*}
whose dimension is the $q$-th persistent sheaf Betti number.
Our goal here is to extend the notion of persistent Laplacian to the setting of cellular sheaves.
We will assume that each stalk of $\mathscr{G}$ is a finite dimensional inner product space, 
and regard any cochain group as a direct sum of inner product spaces. 
Note that the adjoint of $\pi$ is the inclusion map $\iota: C^q(X;\mathscr{F}) \hookrightarrow C^q(Y;\mathscr{G})$. 
We can dualize the above diagram (i.e., reverse all arrows) and define the persistent sheaf Laplacian by the persistent Laplacian of the dualized diagram.
More specifically, we have the following diagram
\[
\begin{tikzcd}
    C^{q-1}(X; \mathscr{F}) \arrow[rr, "d^{q-1}_X", shift left] %\arrow[dd, hook, dashed]
     &
      & C^{q}(X; \mathscr{F}) \arrow[ll, "(d^{q-1}_X)^{\ast}", shift left] \arrow[rd, "d^{q}_{X,Y}", shift left] \arrow[dd, hook, dashed]
       & 
        & \\
     & 
      &
       & \Theta^{q+1}_{X,Y} \arrow[lu, "(d^{q}_{X,Y})^{\ast}", shift left] \arrow[rd, hook, dashed]
        & \\
    %C^{q-1}(Y; \mathscr{G}) \arrow[rr, "d^{q-1}_{Y}"] 
     &
      & C^q(Y; \mathscr{G}) \arrow[rr, "d^q_Y", shift left]
       & 
        & C^{q+1}(Y; \mathscr{G}) \arrow[ll, "(d^q_Y)^{\ast}", shift left]
\end{tikzcd}
\]
(note that an inner product space is self-dual) where $\Theta_{X,Y}^{q+1} = \{x \in C^{q+1}(Y; \mathscr{G}) \mid (d^{q}_{Y})^{\ast}(x) \in C^{q}(X; \mathscr{F})\}$ and 
$d^{q}_{X,Y}$ is the adjoint of $\pi (d^{q}_{Y})^{\ast}\vert_{\Theta^{q+1}_{X,Y}}: \Theta^{q+1}_{X,Y} \to C^q(X;\mathscr{F})$.
We define the $q$-th persistent sheaf Laplacian $\Delta_q^{X,Y}$ by
\begin{align*}
    \Delta_q^{X,Y} = (d^{q}_{X,Y})^{\ast} d^{q}_{X,Y} +  d^{q-1}_X (d^{q-1}_X)^{\ast}.
\end{align*}
The nullity of $\Delta_q^{X,Y}$ is equal to the $q$-th persistent Betti number of the dualized diagram (since we dualize everything, the two cochain complexes become chain complexes).
By the universal coefficient theorem for cohomology, $\pi^{\bullet}$ and $\iota^{\bullet}$ have the same rank (here $\bullet$ means the induced map between cohomology or homology groups).
So the $q$-th persistent Betti number of the dualized diagram is equal to the $q$-th persistent sheaf Betti number.
In other words, we have 
\begin{align*}
    \ker \Delta_q^{X,Y} \cong \pi^{\bullet}(H^q(Y;\mathscr{G})).
\end{align*}
The matrix representation of a persistent sheaf Laplacian can be calculated in the same way as a persistent Laplacian.
\begin{proposition}
    [The matrix representation of a persistent sheaf Laplacian] Choose a basis $\{v_1, \dots, v_n\}$ of
    $\Theta^{q+1}_{X,Y}$ and denote its inner product matrix by $P$. 
    Denote the matrix representation of $(d^q_{X,Y})^{\ast}$ 
    with respect to the canonical basis of $C^q(X; \mathscr{F})$ and $\{v_1, \dots, v_n\}$ as $D^{\ast,q}_{X,Y}$
    and the matrix representation of $d^{q-1}_{X}$ with respect to the canonical bases of $C^{q-1}(X; \mathscr{F})$ 
    and $C^q(X; \mathscr{F})$ as $D^{q-1}_X$. 
    Then, the matrix representation of $\Delta_{q}^{X,Y}$ is 
    \begin{align*}
        D^{q-1}_{X}(D^{q-1}_X)^T + D^{\ast,q}_{X,Y} P^{-1}(D^{\ast,q}_{X,Y})^T.
    \end{align*}
\end{proposition}

\begin{proof}
    We only have to determine the matrix representation $X$ of $d^q_{X,Y}$. 
    Take two vectors $v\in \Theta_{X,Y}^{q+1}, w \in C^q(X; \mathscr{F})$. 
    We abuse the notation a bit and use $v,w$ to denote their coordinates in the form of column vector as well. We have 
    \begin{align}
        \lgl (d^{q}_{X,Y})^{\ast} v,  w\rgl_{C^q(X; \mathscr{F})} &= \lgl v, d^{q}_{X, Y}w \rgl_{\Theta_{X,Y}^{q+1}} \\
        (D_{X,Y}^{\ast, q}v)^Tw &= v^TPXw \\
        v^T(D_{X,Y}^{\ast, q})^Tw &= v^TPXw. 
    \end{align}As $v, w$ are arbitrarily taken, we conclude that $X = P^{-1}(D_{X,Y}^{\ast, q})^T$.
\end{proof}

\begin{proposition}
    The spectrum of the $q$-th persistent sheaf Laplacian does not depend on the orientation of $Y$.
\end{proposition}

\begin{proof}
    Fixing a choice of orientation for each simplex of $Y$,
    it suffices to show that the spectrum of $\Delta_q^{X,Y}$ is unchanged if the orientation of one simplex of $Y$ is alternated.
    We first fix some notations. Suppose we change the orientation of a simplex $\sigma$, then every morphism
    defined with respect to this new orientation will have a bar. We also sometimes denote $\Theta^q_{X,Y}$ by $\Theta$.
    We define a linear map $I_{\sigma, -}$ such that $I_{\sigma, -}\vert_{\mathscr{G}(\sigma)} = -I$ and 
    $I_{\sigma, -}\vert_{\mathscr{G}(\sigma')} = I$ if $\sigma' \neq \sigma$.
    The adjoint of $I_{\sigma, -}$ is itself.
    Depending on context, the domain of $I_{\sigma, -}$ will be understood as 
    $C^{q-1}(X; \mathscr{F}), C^q(Y; \mathscr{G})$ or $C^{q+1}(Y; \mathscr{G})$.
    The proof is divided into cases.
    
    Case I. If $\sigma \in C^{q-1}(X; \mathscr{F})$,
    then $\overline{d^q_{Y}} = d^q_{Y}$ and 
    $\overline{d^{q-1}_X} = d^{q-1}_X I_{\sigma, -}$.
    So $\overline{d^{q-1}_X}(\overline{d^{q-1}_{X}})^{\ast}=
    d^{q-1}_X I_{\sigma, -} I_{\sigma, -}(d^{q-1}_X)^{\ast} = d^{q-1}_X(d^{q-1}_X)^{\ast}$.
    
    Case II. If $\sigma \in C^{q}(Y; \mathscr{G})$,
    then $\overline{d^{q-1}_X} = \pi I_{\sigma, -}\vert_{C^q(X; \mathscr{F})} d^{q-1}_X$
    and $\overline{d^q_{Y}} = d^q_{Y} I_{\sigma, -}$.
    So $(\overline{d^q_{Y}})^{\ast} = I_{\sigma, -} (d^q_{Y})^{\ast}$. As 
    \begin{align*}
        (\overline{d^q_{Y}})^{\ast}\Theta^{q}_{X,Y} = I_{\sigma, -} (d^q_{Y})^{\ast}\Theta^{q}_{X,Y} \subset I_{\sigma, -}C^q(X; \mathscr{F}) = C^q(X; \mathscr{F}),
    \end{align*}we see that $\Theta^{q}_{X,Y} \subset \overline{\Theta^q_{X,Y}}$. 
    Similarly  
    \begin{align*}
        (d^q_{Y})^{\ast} \overline{\Theta^q_{X,Y}} = I_{\sigma, -}I_{\sigma, -}(d^q_{Y})^{\ast} \overline{\Theta^q_{X,Y}} 
        = I_{\sigma, -} (\overline{d^q_{Y}})^{\ast}\overline{\Theta^q_{X,Y}} \subset I_{\sigma, -}C^q(X;\mathscr{F})=C^q(X; \mathscr{F}),
    \end{align*}
    we see that $\Theta^{q}_{X,Y} \supset \overline{\Theta^{q}_{X,Y}}$, so $\Theta^{q}_{X,Y} = \overline{\Theta^{q}_{X,Y}}$.

    Then 
    \begin{align*}
        \pi (\overline{d^q_{Y}})^{\ast} \vert_{\Theta^{q}_{X,Y}} = \pi I_{\sigma, -}\vert_{C^q(X; \mathscr{F})} \pi (d^q_{Y})^{\ast} \vert_{\Theta^{q}_{X,Y}}.
    \end{align*}
    So 
    \begin{align*}
        \overline{\Delta_q^{X,Y}} &= \overline{d^{q-1}_X} (\overline{d^{q-1}_X})^{\ast} + \pi (\overline{d^q_{Y}})^{\ast} \vert_{\Theta^{q}_{X,Y}} \big(\pi (\overline{d^q_{Y}})^{\ast} \vert_{\Theta^{q}_{X,Y}}\big)^{\ast} \\
        &= \pi I_{\sigma, -}\vert_{C^q(X;\mathscr{F})} \Delta_q^{X,Y} (\pi I_{\sigma, -}\vert_{C^q(X;\mathscr{F})})^{\ast}.
    \end{align*}
    
    Case III. If $\sigma \in C^{q+1}(Y; \mathscr{G})$,
    then $\overline{d^{q-1}_X} = d^{q-1}_X$ and $\overline{d^q_{Y}} = I_{\sigma, -} d^q_{Y}$.
    So 
    $(\overline{d^q_{Y}})^{\ast} = (d^q_{Y})^{\ast} I_{\sigma, -}$.
    As 
    $C^q(X;\mathscr{F}) \supset (\overline{d^q_{Y}})^{\ast} \overline{\Theta} = (d^q_{Y})^{\ast}
    I_{\sigma, -}\overline{\Theta}$, we see that 
    $I_{\sigma, -}\overline{\Theta} \subset \Theta$. 
    Similarly 
    $C^q(X;\mathscr{F}) \supset (d^q_{Y})^{\ast} \Theta = (\overline{d^q_{Y}})^{\ast} I_{\sigma,-}\Theta$,
    implying that 
    $I_{\sigma, -}\Theta \subset \overline{\Theta}$. 
    So $I_{\sigma, -}\Theta = \overline{\Theta}$.
    Denote by $I_{\sigma,-}^{\overline{\Theta}}: \overline{\Theta} \to \Theta$ the restriction of $I_{\sigma,-}$ on $\overline{\Theta}$.
    We have $\pi (\overline{d^q_{Y}})^{\ast}\vert_{\overline{\Theta}} = \pi (d_Y^q)^{\ast}I_{\sigma,-}\vert_{\overline{\Theta}} = \pi(d_Y^q)^{\ast}\vert_{\Theta}I_{\sigma,-}^{\overline{\Theta}}$. Then
    \begin{align*}
        \overline{\Delta_q^{X,Y}} &= \overline{d^{q-1}_X} (\overline{d^{q-1}_X})^{\ast} + \pi(\overline{d^q_{Y}})^{\ast}\vert_{\overline{\Theta}} \big(\pi (\overline{d^q_{Y}})^{\ast} \vert_{\overline{\Theta}}\big)^{\ast} \\
        &= d^{q-1}_X (d^{q-1}_X)^{\ast} + \pi(d^q_{Y})^{\ast}\vert_{\Theta} I_{\sigma, -}^{\overline{\Theta}} (I_{\sigma, -}^{\overline{\Theta}})^{\ast}
        \big(\pi (d^q_{Y})^{\ast}\vert_{\Theta}\big)^{\ast}\\
        &= d^{q-1}_X (d^{q-1}_X)^{\ast} + \pi(d^q_{Y})^{\ast}\vert_{\Theta}
        \big(\pi (d^q_{Y})^{\ast}\vert_{\Theta}\big)^{\ast}\\
        &= \Delta_q^{X,Y}.
    \end{align*}
\end{proof}

\begin{example}
    Consider the 1-dimensional simplicial complex $Y$
    \begin{align*}
        \begin{tikzpicture}
            \centering
            %% vertices
            \draw[fill=black] (0,0) circle (1pt);
            \draw[fill=black] (6,0) circle (1pt);
            \draw[fill=black] (2,0) circle (1pt);
            \draw[fill=black] (4,0) circle (1pt);
            %% vertex labels
            \node at (0,-0.5) {1};
            \node at (6,-0.5) {2};
            \node at (2,-0.5) {3};
            \node at (4,-0.5) {4};
            %%% edges
            \draw[thick] (0,0) -- (2,0) -- (4,0) -- (6,0);
            %\caption{}
        \end{tikzpicture}
    \end{align*}
    and the constant sheaf $\underline{\mathbb{R}}$ over $Y$. 
    We compute $\Delta_{0}^{X,Y}$ when $X= \{1,2\}$. 
    The matrix representation of $(d^0_{Y})^{\ast}$ 
    \begin{align*}
        \bordermatrix{
            ~ & 13 & 34 & 42  \cr
            1 & -1 &  0 &  0  \cr
            2 &  0 &  0 &  1  \cr
            3 &  1 & -1 &  0  \cr
            4 &  0 &  1 & -1  \cr
        }.
    \end{align*} 
    After a few steps of Gauss elimination we get
    \begin{align*}
        (d^0_Y)^{\ast} = \bordermatrix{
            ~ & 13 & 13+34 & 13+34+42  \cr
            1 & -1 &    -1 &       -1  \cr
            2 &  0 &     0 &        1  \cr
            3 &  1 &     0 &        0  \cr
            4 &  0 &     1 &        0  \cr
        }.
    \end{align*}
    It is clear that $\Theta^{0}_{X,Y}=\text{span}\{13+34+42\}$, $P=3$ and 
    the matrix representation of $(d^{0}_{X,Y})^{\ast}$ is
    \begin{align*}
        \bordermatrix{
            ~ & 13+34+42  \cr
            1 &       -1  \cr
            2 &        1  \cr
        }.
    \end{align*} Then, the matrix representation of $\Delta_{0}^{X,Y}$ is 
    \begin{align*}
        \begin{pmatrix}
            1/3 & -1/3 \\
            -1/3 & 1/3
        \end{pmatrix}
    \end{align*} and its spectrum is $\{0,2/3\}$. 
\end{example}

\section{Experiments} \label{experiments}

We first utilize some simple shapes to illustrate persistent sheaf Laplacians. Applications to realistic molecules are also given. 

\subsection{Simple shapes}

\begin{figure}[htbp] 
    \centering
    \begin{subfigure}{0.5\textwidth}
        \centering
        \includegraphics[width=0.5\textwidth]{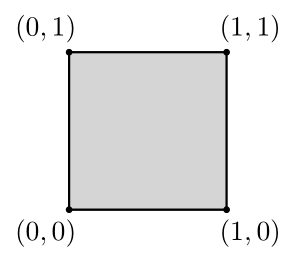}
        %\caption{Pentagon}
    \end{subfigure}\hfill
    \begin{subfigure}{0.5\textwidth}
        \centering
        \includegraphics[width=0.58\textwidth]{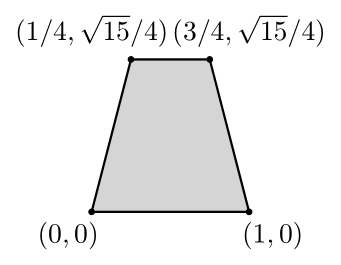} 
        %\caption{Heptagon}
    \end{subfigure}\hfill
    \caption{A square and a trapezoid. Coordinates of vertices are shown.}
    \label{2d-structures}
\end{figure}

Given a labeled point cloud $P$
(i.e., a point cloud with a nonzero quantity $q_i$ associated with each point $v_i$),
we can build a Rips or alpha filtration out of it and construct a sheaf $\mathscr{S}_t$ on each $X_t$ consistently
as described in section \ref{cellular sheaves on a labeled simplicial complex} provided a suitable global $F: 2^P \to \real$ is chosen. 
In other words, for a face relation $\sigma \leqslant \tau$, the stalks of $\sigma$, $\tau$ and the restriction map 
remain the same for any $\mathscr{S}_t$ containing them. 
Since $\mathscr{S}_t$ is the pullback of $\mathscr{S}_{t+p}$ for any $t$ and $p$,  
we get a persistent module of sheaf cochain complexes and can compute the spectra of persistent sheaf Laplacians. 
In this section, we calculate the spectra of persistent sheaf Laplacians for a few examples of point clouds in this way. 
Some examples are the vertices of simple geometrical shapes and some are the coordinates of atoms of molecules.
We assign the quantities $q_i$ to simple geometrical shapes, and take the partial charges as $q_i$ for molecules.
An alpha filtration is built for each labeled point cloud, parametrized by radius $r$.  
We choose $F$ such that $F$ maps every vertex $v_i$ to 1, every edge $v_iv_j$ to the length of itself $r_{ij}$, and every 2-cell $v_iv_jv_k$ to the product of lengths of its edges $r_{ij}r_{ik}r_{jk}$. 
The spectra of $\Delta_q^{X_r, X_{r+p}}$ for $q=0,1$ and some selected $r,p$ will be calculated.
The radius $r$ will be a multiple of 0.01 or 0.01\AA.
The minimal nonzero eigenvalue of the persistent sheaf Laplacian $\Delta_q^{X_r,X_{r+p}}$ is denoted by $\lambda_q^{r, p}$
and the $q$-th persistent sheaf Betti number of the pair $X_r$, $X_{r+p}$ is denoted by $\beta_q^{r, p}$.
The first set of examples are a 2-dimensional square and a 2-dimensional trapezoid as shown in Figures \ref{2d-structures} with different choices of local property $q_i$.

\begin{figure}[htbp]
    \centering
    \begin{subfigure}[b]{0.9\textwidth}
        \centering
        \includegraphics[width=0.9\linewidth]{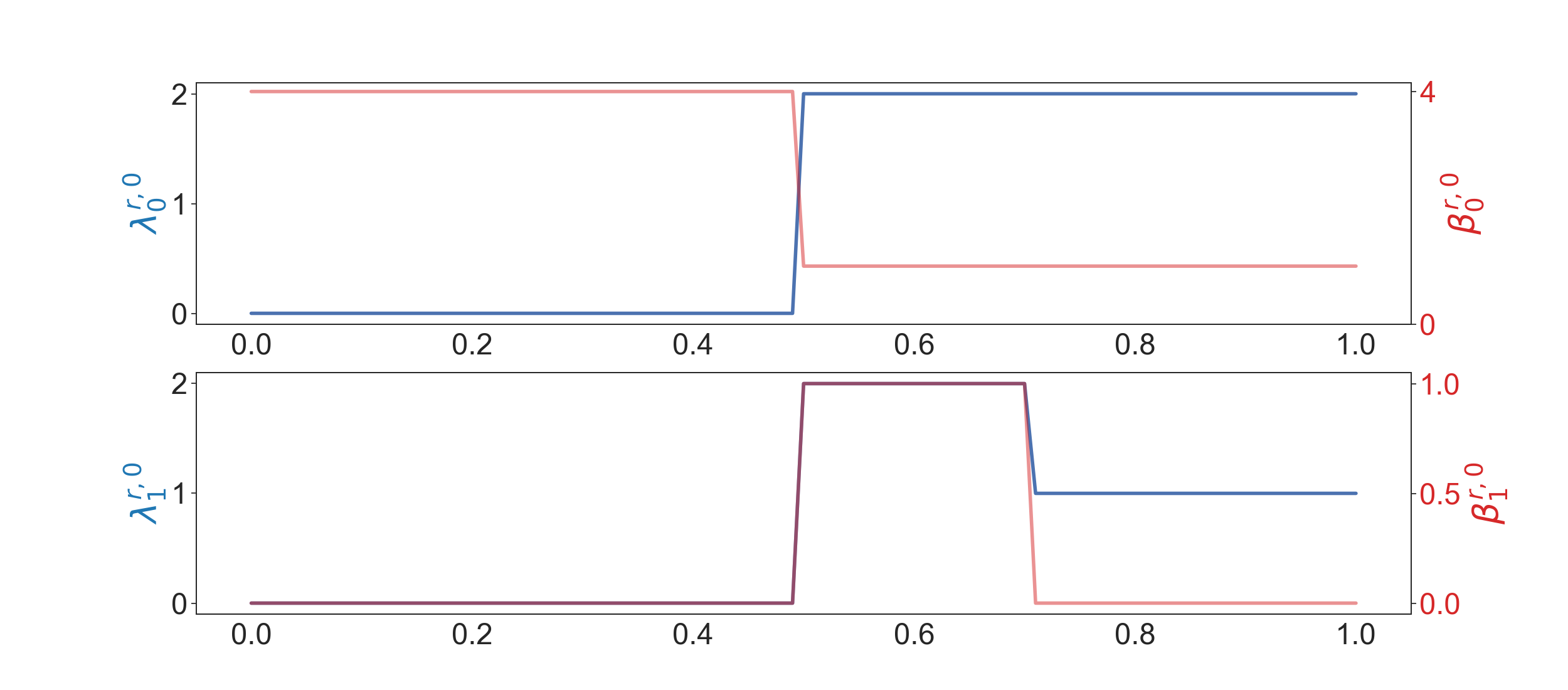}
    \end{subfigure}
    \begin{subfigure}[b]{0.9\textwidth}
        \centering
        \includegraphics[width=0.9\linewidth]{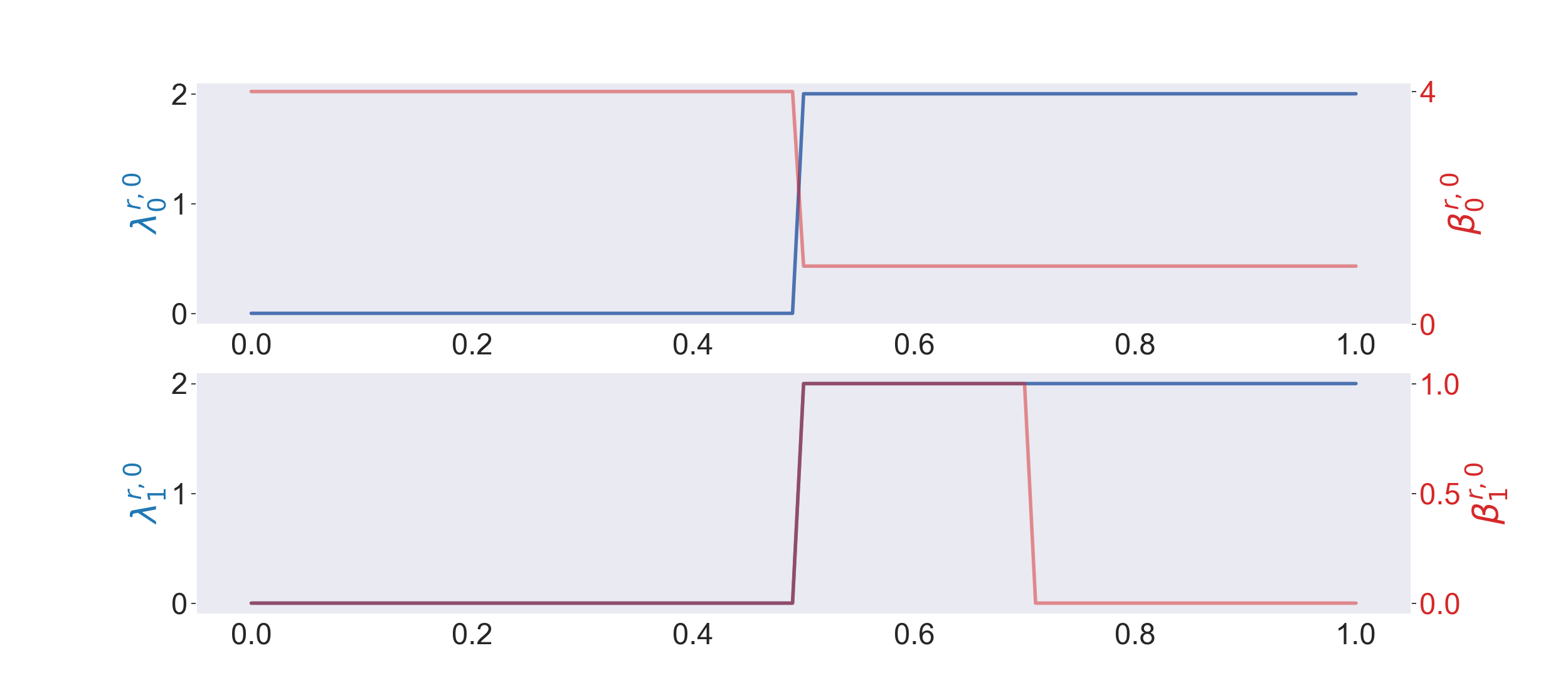}
    \end{subfigure}
    \caption{The results of the square with charge 1 assigned to each vertex when $p=0$ are shown in figures with grey background. 
    The results of the square if each $\mathscr{S}_t$ is a constant sheaf $\underline{\mathbb R}$ when $p=0$ are shown in figures with white background. 
    Here $\beta_{q}^{r, 0}, \lambda_{q}^{r, 0}$ are the persistent sheaf Betti number and the minimal nonzero eigenvalues of $\Delta_q^{X_r, X_r}$.   
    }
    \label{square_p=0_++++}
\end{figure}

\begin{figure}[htbp]
    \centering
    \begin{subfigure}[b]{0.9\textwidth}
        \centering
        \includegraphics[width=0.9\linewidth]{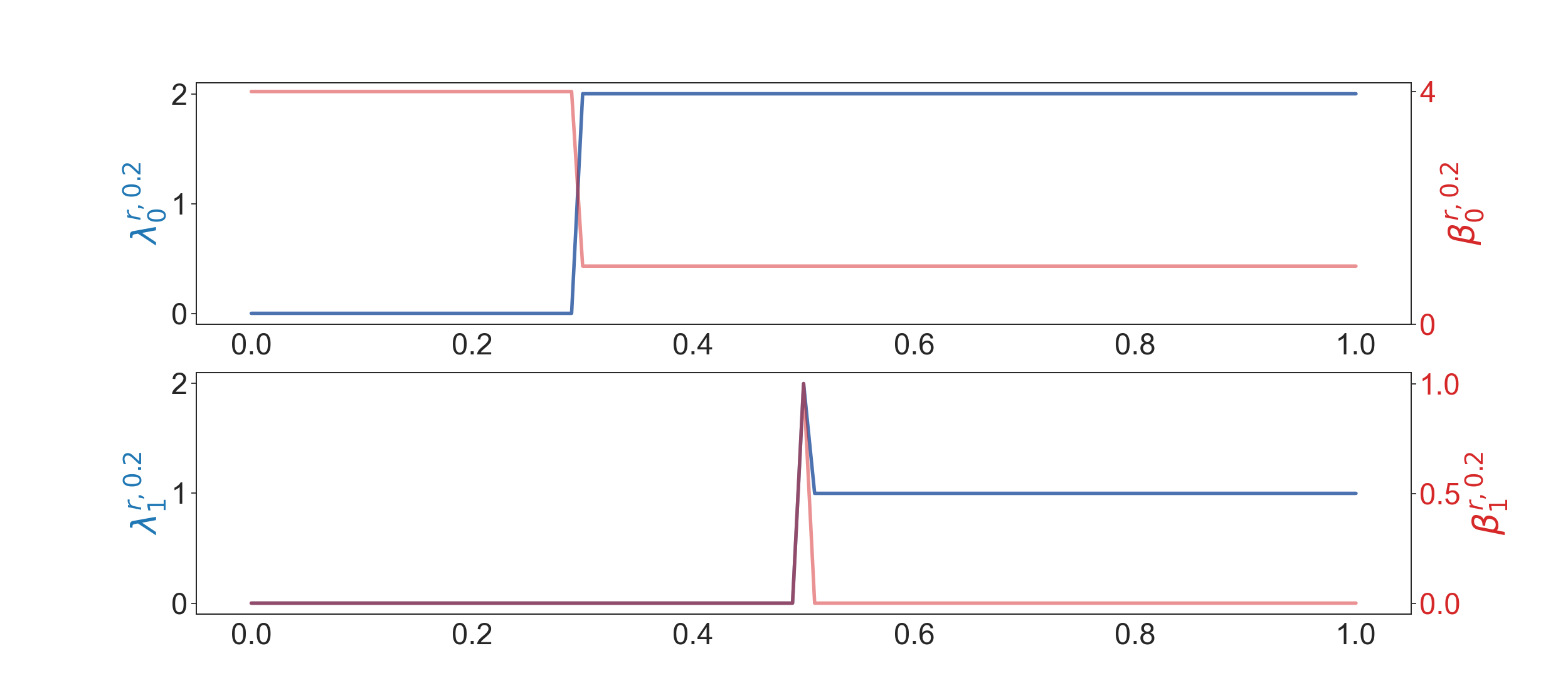}
    \end{subfigure}
    \begin{subfigure}[b]{0.9\textwidth}
        \centering
        \includegraphics[width=0.9\linewidth]{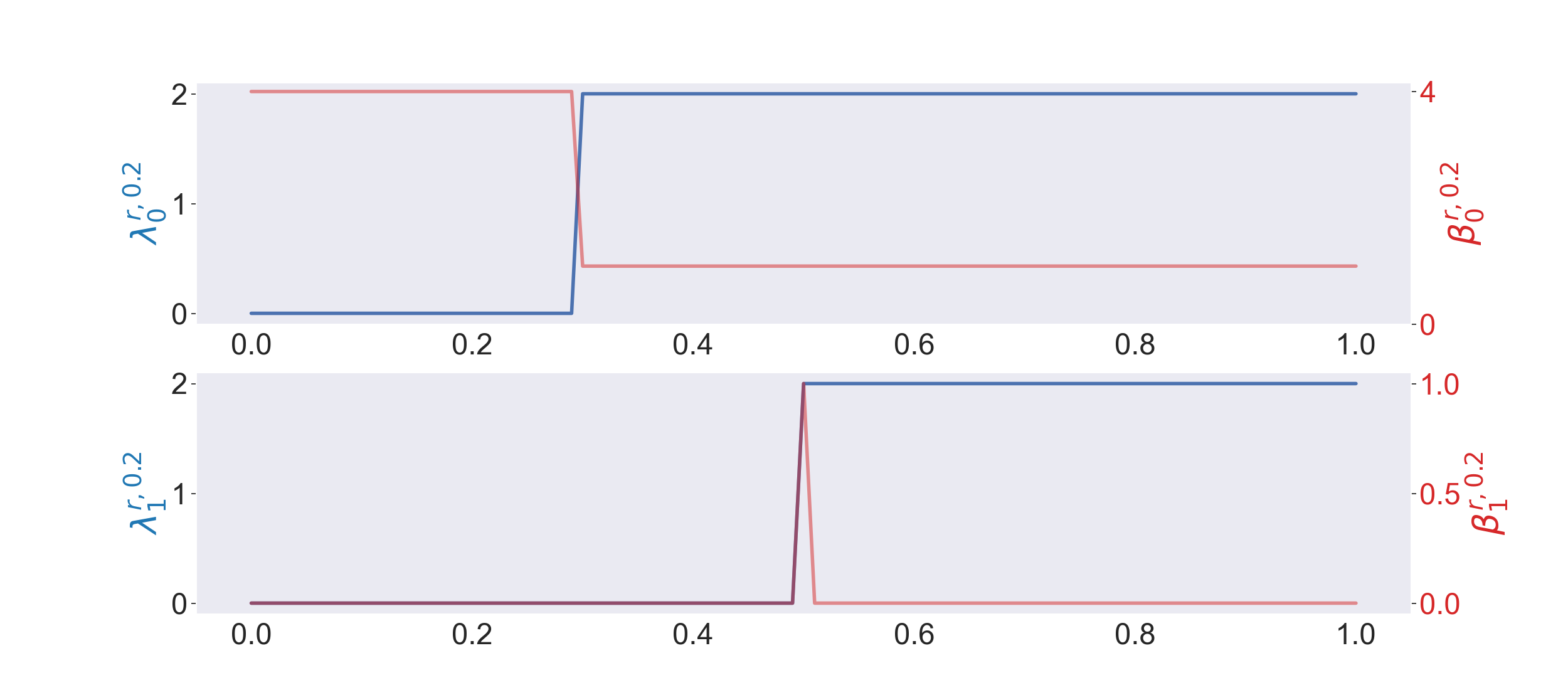}
    \end{subfigure}
    \caption{The results of the square with charge 1 assigned to each vertex when $p=0.2$ are shown in figures with grey background. 
    The results of the square if each $\mathscr{S}_t$ is a constant sheaf $\underline{\mathbb R}$ when $p=0$ are shown in figures with white background.    }
    \label{square_p=0.2_++++}
\end{figure}

\begin{figure}[htbp]
    \centering
    \begin{subfigure}[b]{0.9\textwidth}
        \centering
        \includegraphics[width=0.9\linewidth]{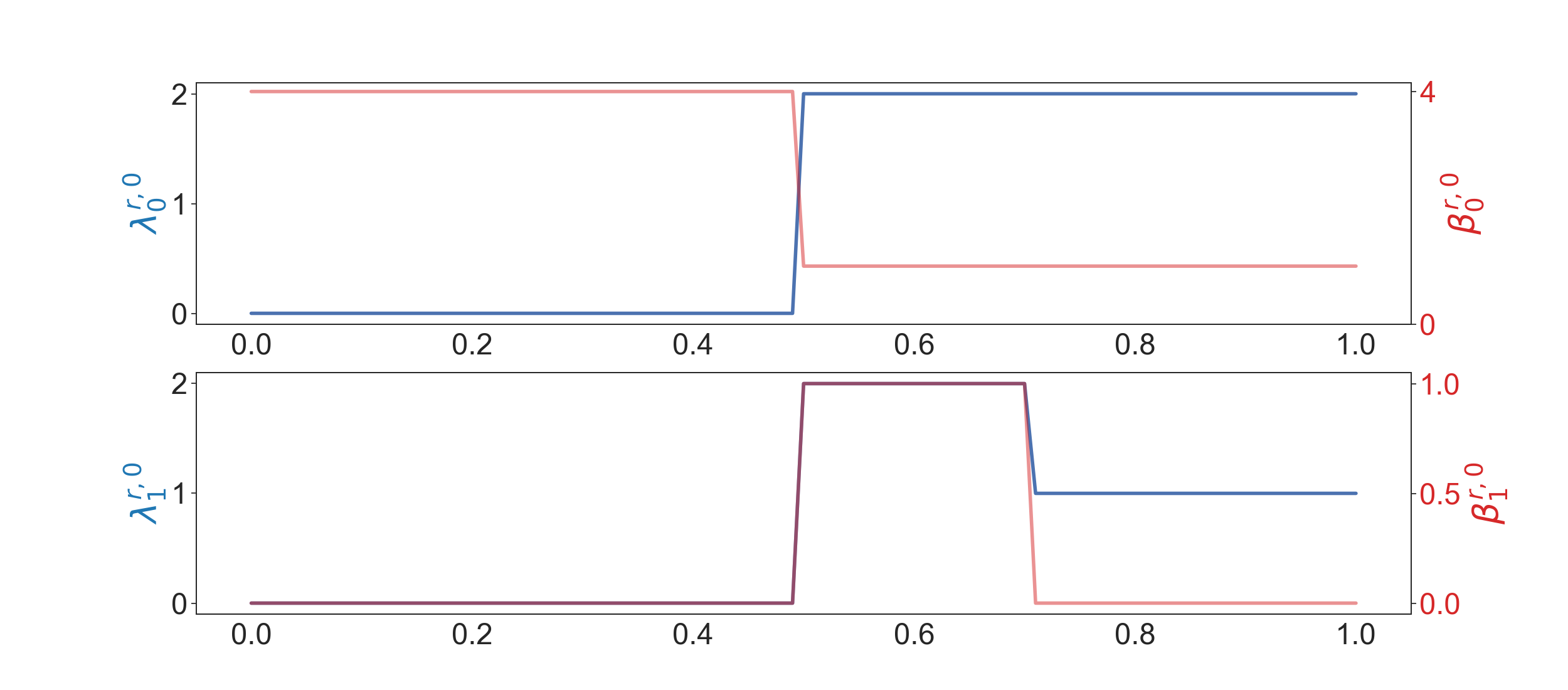}
    \end{subfigure}
    \begin{subfigure}[b]{0.9\textwidth}
        \centering
        \includegraphics[width=0.9\linewidth]{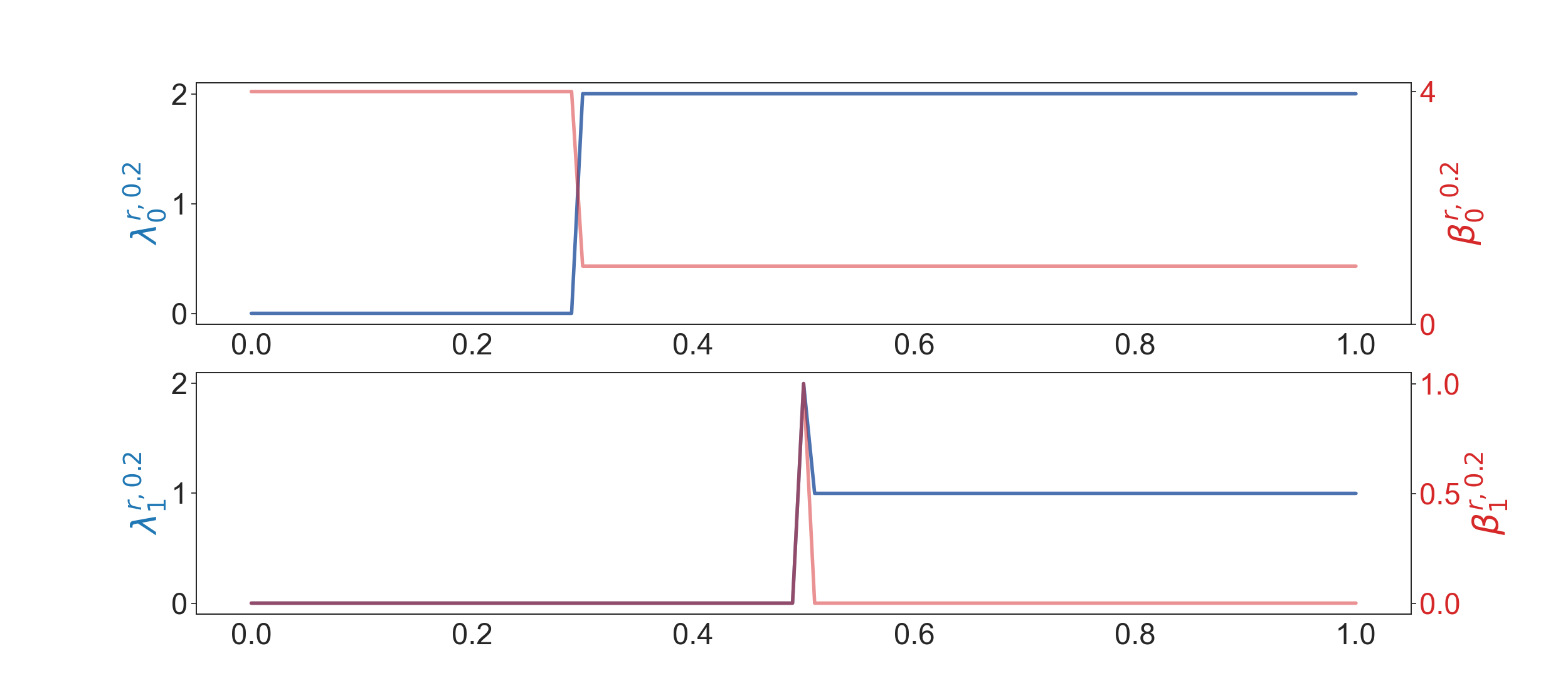}
    \end{subfigure}
    \caption{The results of the square with charge 1 assigned to $(0,0), (1,0), (1,1)$
    and -1 assigned to $(0,1)$ when $p=0$ and $p=0.2$. 
    }
    \label{square_+++-}
\end{figure}

\begin{figure}[htbp]
    \centering
    \begin{subfigure}[b]{0.9\textwidth}
        \centering
        \includegraphics[width=0.9\linewidth]{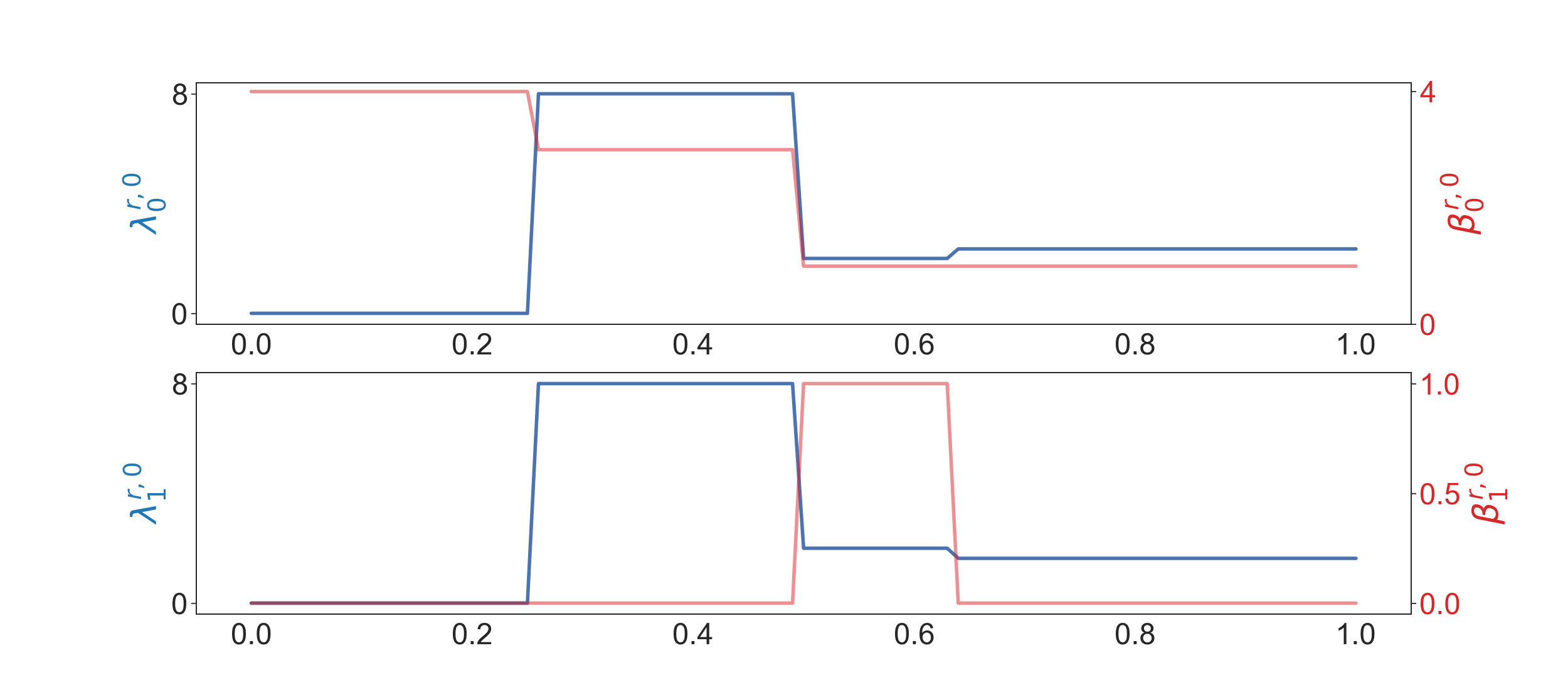}
    \end{subfigure}
    \begin{subfigure}[b]{0.9\textwidth}
        \centering
        \includegraphics[width=0.9\linewidth]{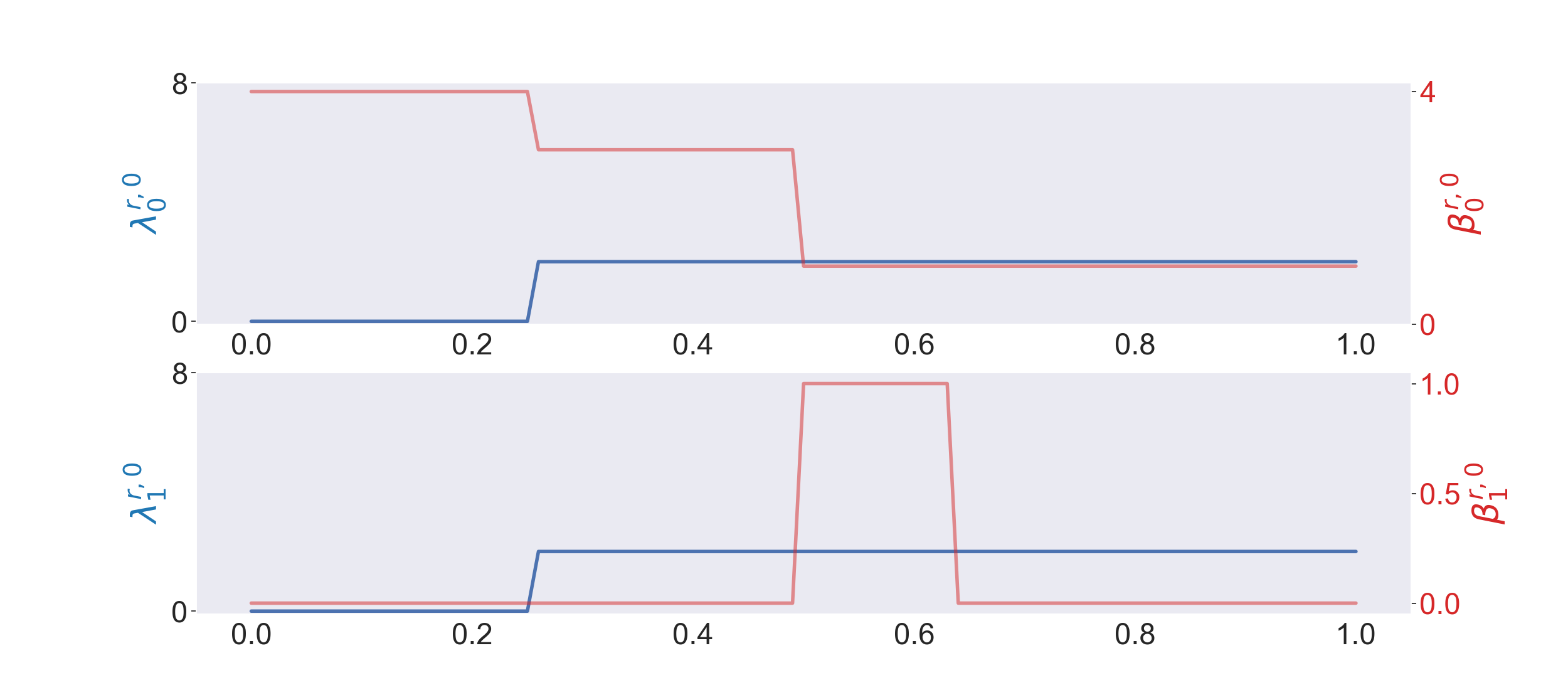}
    \end{subfigure}
    \caption{The results of the trapezoid with charge 1 assigned to each vertex when $p=0$ are shown in figures with grey background. 
    The results of the trapezoid with the constant sheaf $\underline{\mathbb R}$ when $p=0$ are shown in figures with white background. 
    }
    \label{trapezoid_p=0_++++}
\end{figure}

\begin{figure}[htbp]
    \centering
    \begin{subfigure}[b]{0.9\textwidth}
        \centering
        \includegraphics[width=0.9\linewidth]{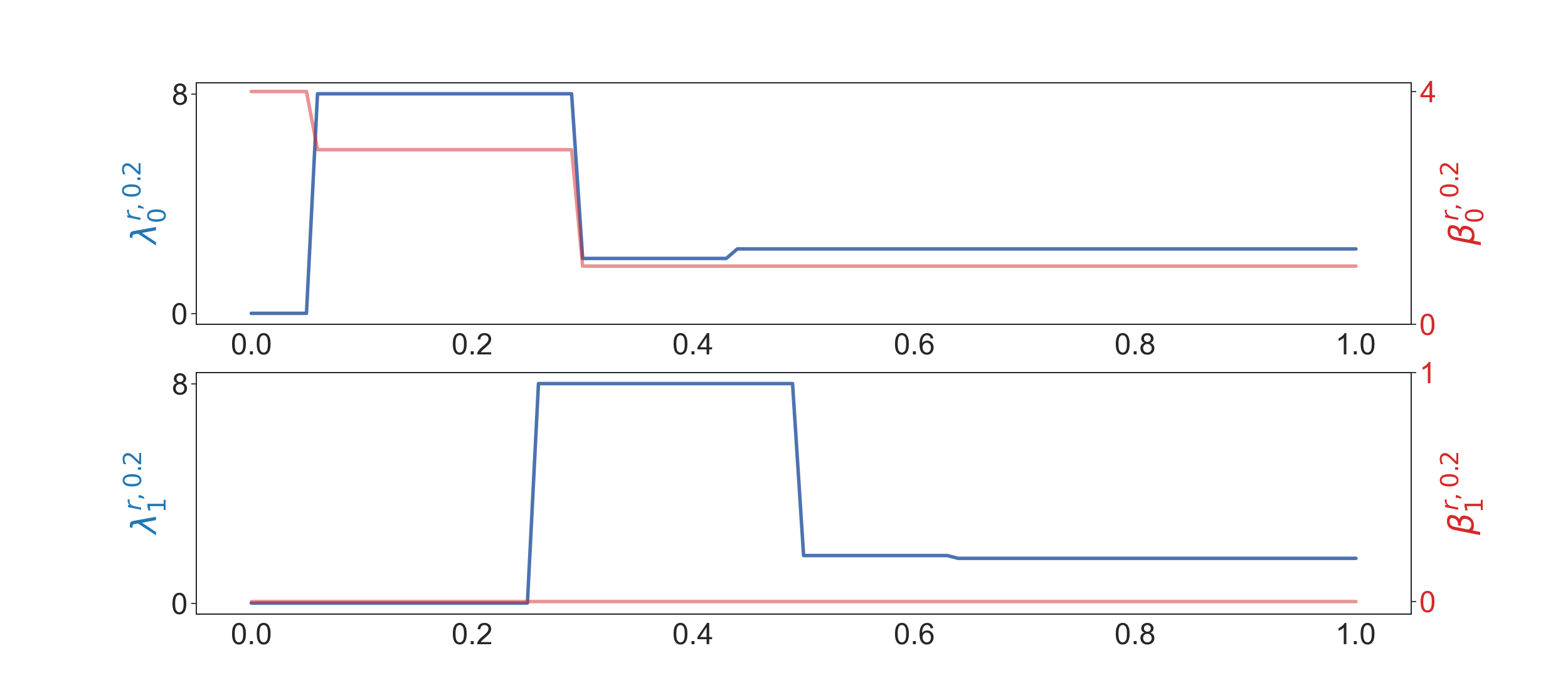}
    \end{subfigure}
    \begin{subfigure}[b]{0.9\textwidth}
        \centering
        \includegraphics[width=0.9\linewidth]{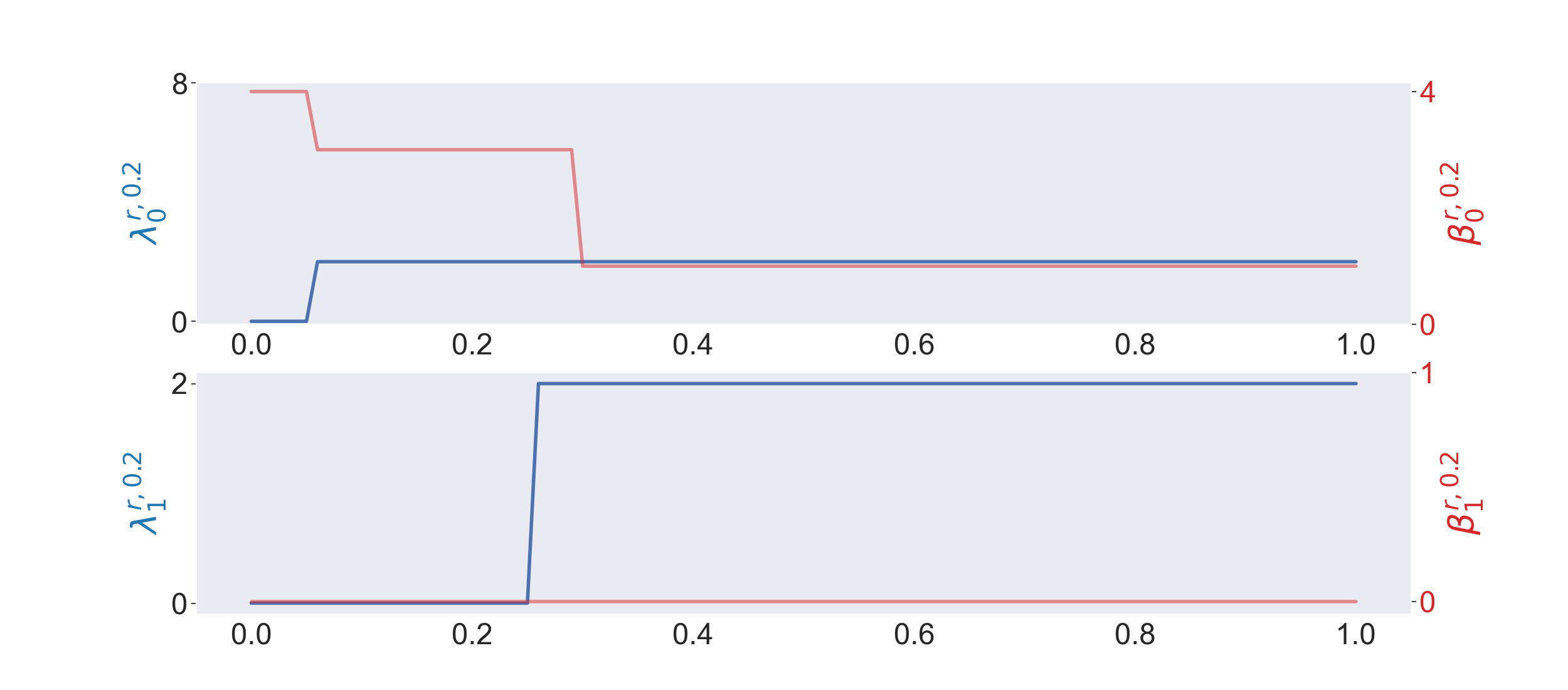}
    \end{subfigure}
    \caption{The results of the trapezoid with charge 1 assigned to each vertex when $p=0.2$ are shown in figures with white background. 
    The results of the trapezoid with the constant sheaf $\underline{\mathbb R}$ when $p=0.2$ are shown in figures with grey background.}
    \label{trapezoid_p=0.2_++++}
\end{figure}

\begin{figure}[htbp]
    \centering
    \includegraphics[width=0.9\linewidth]{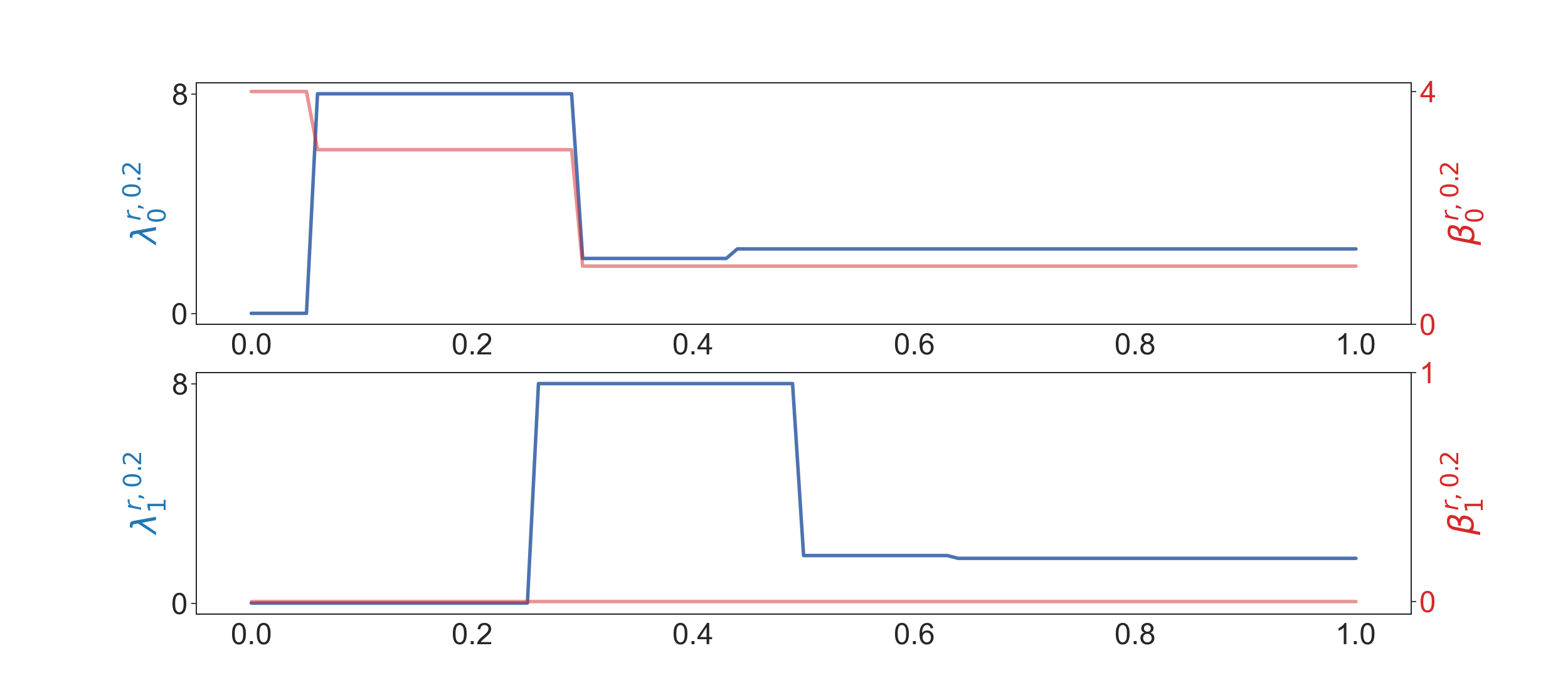}
    \caption{The results of the trapezoid with charge 1 assigned to $(0,0), (1,0), (3/4, \sqrt{15}/4)$ and -1 
    assigned to $(1/4, \sqrt{15}/4)$ when $p=0.2$. 
    }
    \label{trapezoid_p=0.2_+++-}
\end{figure}

The results of the square and the trapezoid are shown in Figures 
\ref{square_p=0_++++},
\ref{square_p=0.2_++++}, 
\ref{square_+++-}, 
\ref{trapezoid_p=0_++++}, 
\ref{trapezoid_p=0.2_++++} and 
\ref{trapezoid_p=0.2_+++-}. 
For a comparison, the results when each sheaf in the filtration is a constant sheaf are also presented. 

Figure \ref{square_+++-} gives a case where one of the four nodes flips its charge from 1 to -1.  
We observe that the change of signs of $q_i$ does not affect the results, 
though the eigenspaces of Laplacians will be different.

For the trapezoid structures, we see patterns very similar those of squares expect for additional complexity reflecting the geometry.
In both types of geometries, there are certain  similarity and difference between the spectra of sheaves and constant sheaves, indicating these methods  may complement  in practical applications. 

\begin{figure}[htbp]
    \centering
    \includegraphics[width=0.3\linewidth]{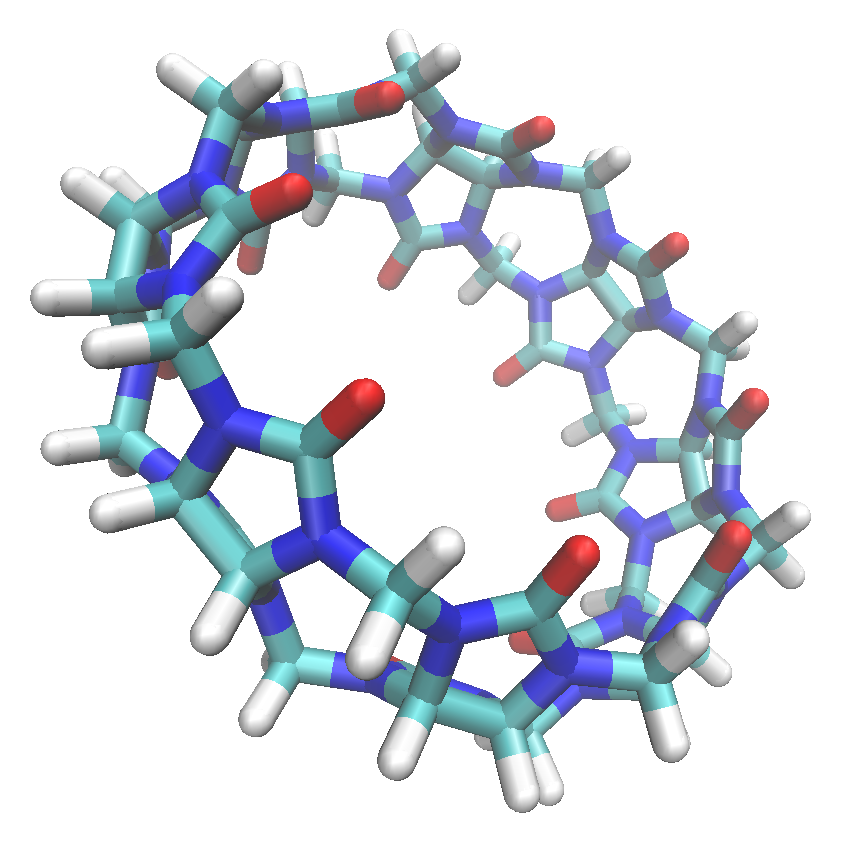}
    \caption{The structure of cucurbit[8]uril. 
    }
    \label{cucurbit}
\end{figure}

\begin{figure}[htbp]
    \centering
    \includegraphics[width=0.9\linewidth]{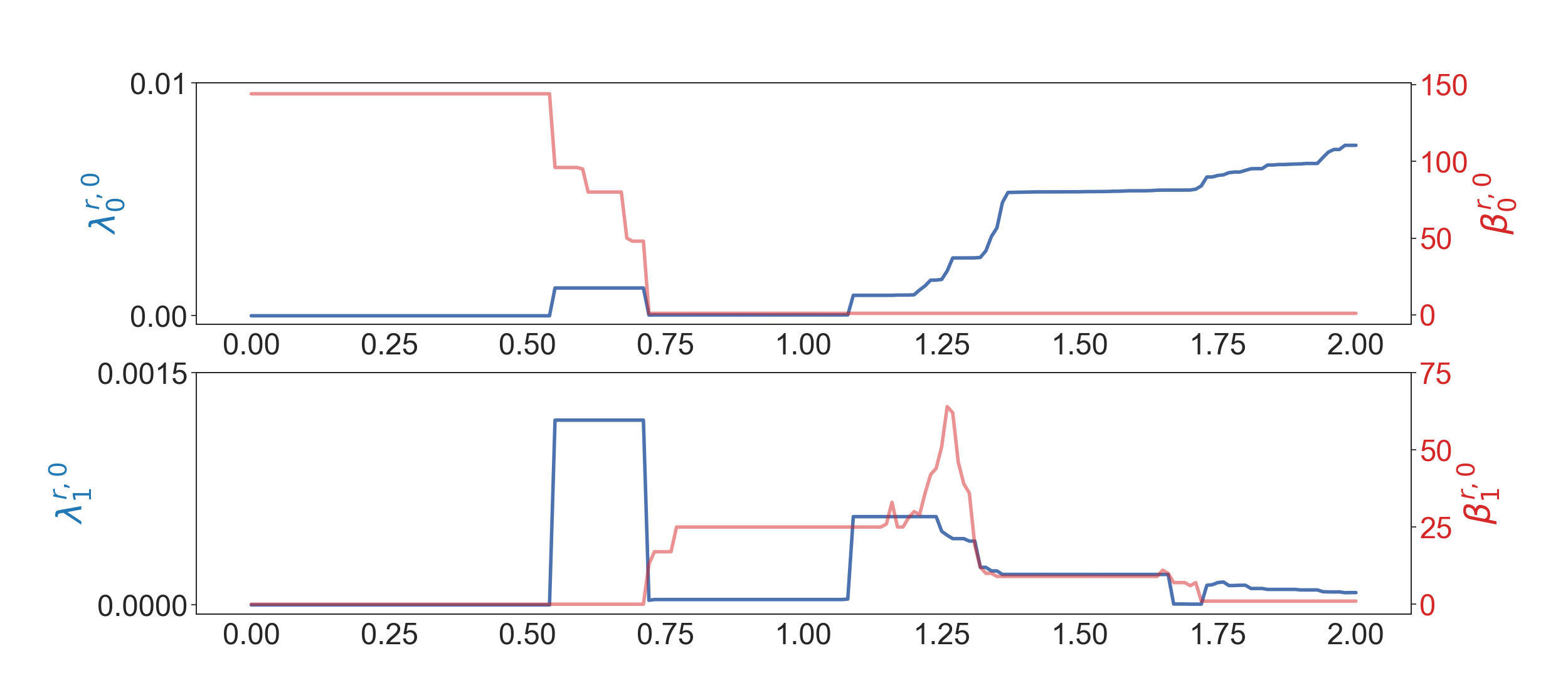}
    \caption{The results of the cucurbit[8]uril when $p=0$. 
    }
    \label{cb8_p=0}
\end{figure}

\begin{figure}[htbp]
    \centering
    \includegraphics[width=0.9\linewidth]{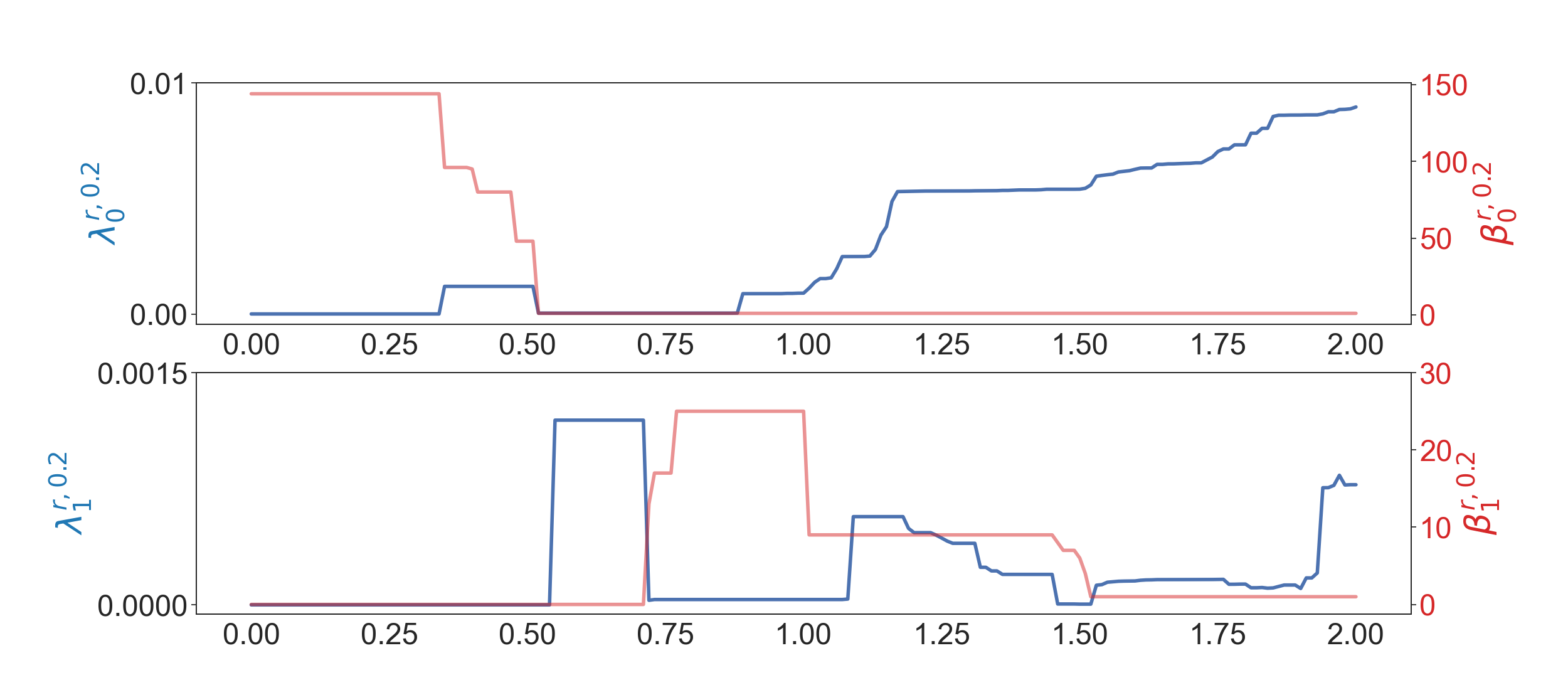}
    \caption{The results of the cucurbit[8]uril when $p=0.2$. 
    }
    \label{cb8_p=0.2}
\end{figure}

\subsection{Complex molecules}

Next we study the molecule cucurbit[8]uril \cite{samplchallenges}, which is shown in Figure \ref{cucurbit}.
We associate each atom with the corresponding partial charge (obtained using \cite{ravcek2020atomic}).
The results for cucurbit[8]uril are shown in Figures \ref{cb8_p=0} and \ref{cb8_p=0.2}. Due to complexity of the molecule, it is very difficult to explain the spectral details of the system. 
However, this information can be very useful for machine learning analysis. 

\begin{figure}[htbp]
    \centering
    \includegraphics[width=0.3\linewidth]{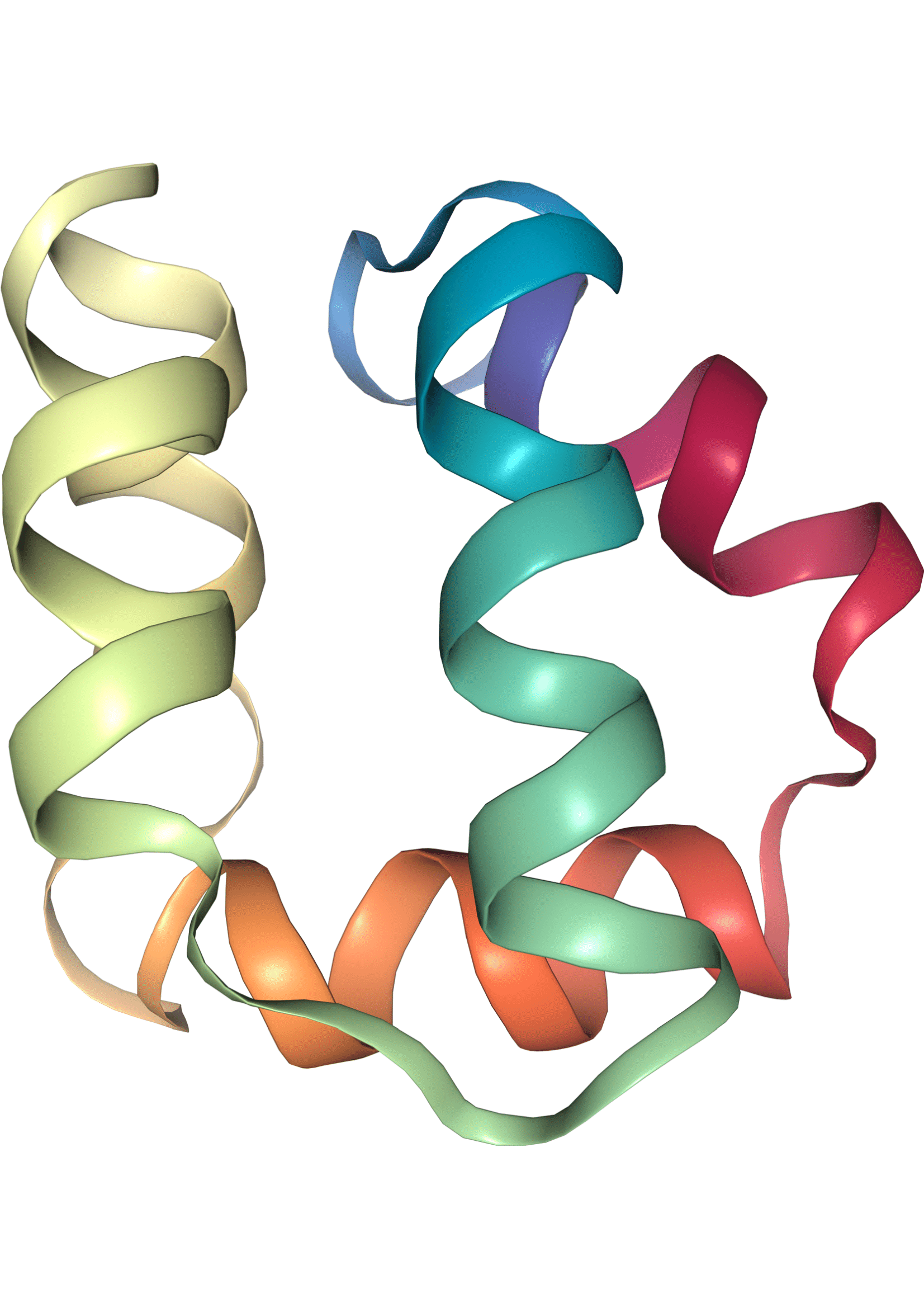}
    \caption{Illustration of the structure of bacteriocin AS-48. 
    }
    \label{1e68}
\end{figure}

\begin{figure}[htbp]
    \centering
    \includegraphics[width=0.9\linewidth]{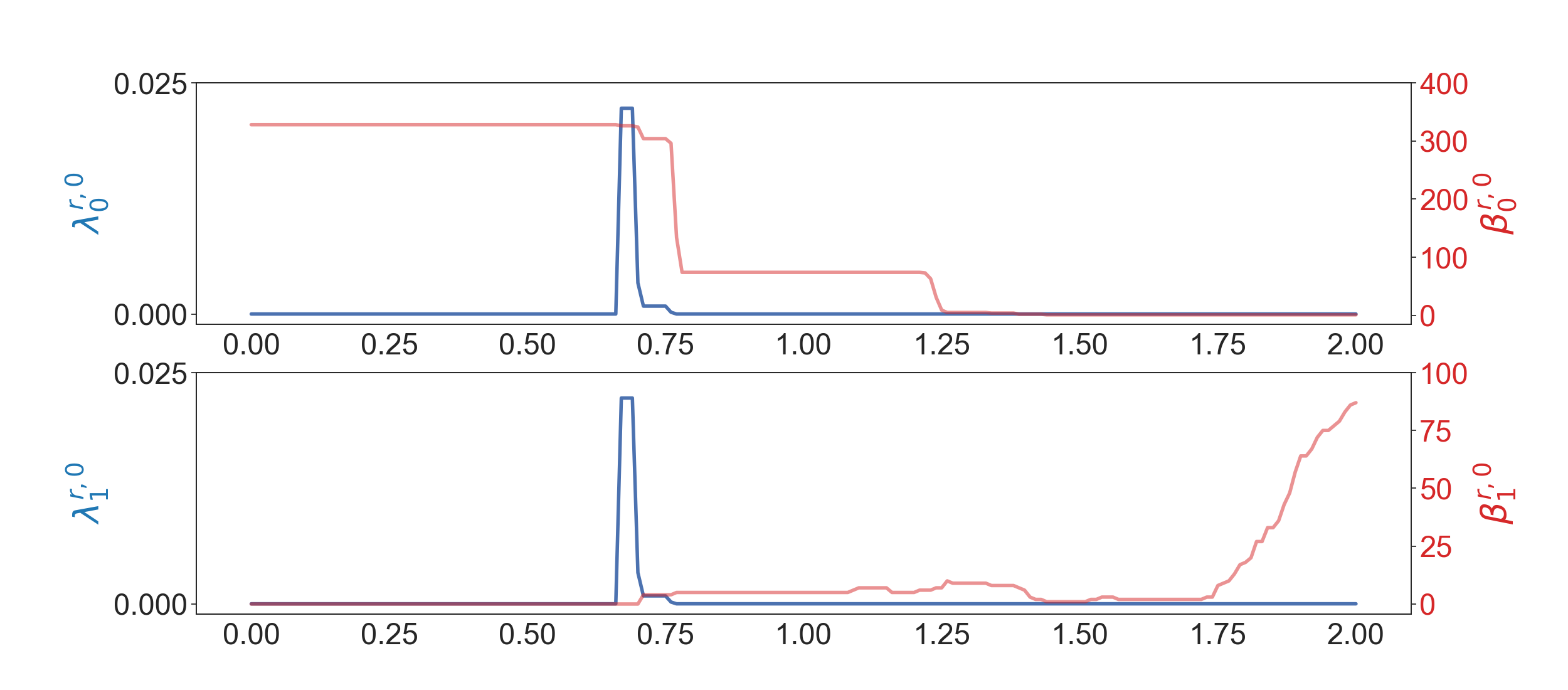}
    \caption{The results of AS-48 when $p=0$. 
    }
    \label{1e68_p=0}
\end{figure}

\begin{figure}[htbp]
    \centering
    \includegraphics[width=0.9\linewidth]{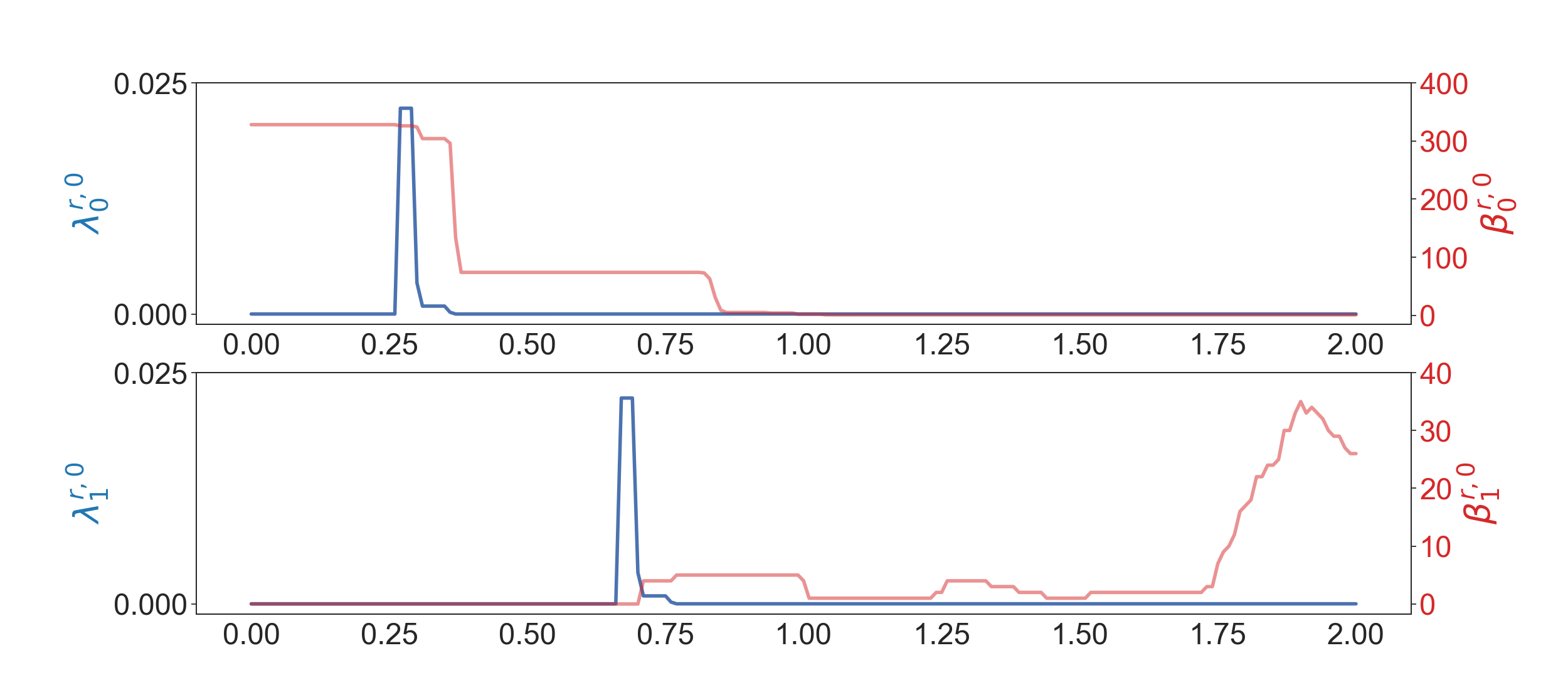}
    \caption{Results of AS-48 when $p=0.4$.
    }
    \label{1e68_p=0.4}
\end{figure}

Finally, to demonstrate our method for practical problems, we study a small protein called bacteriocin AS-48 (PDB ID: 1E68) \cite{1e68}. 
We select the model 1  of AS-48 and compute the pqr file by PDB2PQR with the Amber force field \cite{Jurrus2018}.
For the sake of faster computation, we only use the coordinates of carbon atoms as the point cloud.
Results are shown in Figures \ref{1e68_p=0} and \ref{1e68_p=0.4}.

\section{Concluding remarks}

The success of persistent homology in data science lies in its integration of multiscale analysis and topological abstraction in terms of topological invariants \cite{edelsbrunner2008persistent,zomorodian2005computing}.  
However, this approach has a significant drawback, i.e., its inability to detect homotopic shape evolution of data that involves  no topological changes.  
This limitation was addressed by persistent topological Laplacians, including persistent spectral graph, also called persistent Laplacians \cite{wang2020persistent}, and evolutionary de Rahm Hodge theory, also called persistent Hodge Laplacians \cite{chen2019evolutionary}.  
The multiscale geometric analysis in these new approaches has dramatically empowered classical spectral theory. 
The harmonic spectra of persistent Laplacians fully retain the topological invariants of persistent homology, 
while their non-harmonic spectra capture the homotopic shape evolution of data that cannot be obtained from persistent homology.  
Nonetheless, persistent Laplacians cannot handle non-geometric information of the data, such as the atomic property of a molecule and node property of a network.  
This work proposes persistent sheaf Laplacians (PSL) to overcome this limitation. 
The theory of PSLs introduces multiscale analysis to cellular sheaf Laplacians to further extend their application potential via a filtration. 
{ The proposed PSLs} have been demonstrated on various point cloud data with localized labels.  
%It addresses the need to encode non-geometrical information, such as element specific multiscale information \cite{cang2018integration}.
%Stated differently, persistent sheaf Laplacians  graft   persistent spectral graph  with a sheaf structure, which substantially empowers   persistent   Laplacians. 
%Practical examples are discussed for the Schr\"{o}dinger equation, Coulomb’s potential, and a protein point cloud data to illustrate the proposed mathematical tools. 

Potential future development is widely open. 
First, a fast implementation of PSLs and a good understanding of the eigenvalues and eigenspaces of PSLs will dramatically boost the applications of PSLs to complex and big data.
Second, the incorporation of foundational laws of physics, principles of chemistry, and rules of life into PSLs to study various types of data will be exciting. 
Third, the application of PSLs for data fusion, i.e., integrating multiple data sources to produce more consistent, accurate, and useful information than that concatenates individual data sources, is another direction. 
Fourth, we need to explore various specific types of sheaf Laplacians for each specific problem and test their power in applications. 
Fifth, a persistent sheaf Dirac theory will extend the Dirac formulation of quantum persistent homology \cite{ameneyro2022quantum} and provide a potentially more powerful approach for data science.  
Finally, in the spirit of evolutionary de Rahm-Hodge theory \cite{chen2019evolutionary}, one can develop  evolutionary sheaf Dirac on manifolds for volumetric data. 
However, like persistent Hodge Laplacians, this approach can be computationally demanding to implement.

\section*{Acknowledgments}
This work was supported in part by NIH grants R01GM126189, R01AI164266, and R35GM148196, National Science Foundation grants DMS2052983, DMS-1761320, and IIS-1900473, NASA  grant 80NSSC21M0023,   Michigan State University Research Foundation, and  Bristol-Myers Squibb  65109.  

\clearpage 
%\printbibliography
\bibliographystyle{abbrv}
\bibliography{xiaoqi}

\end{document}